\documentclass{eajam}
\usepackage{charter}
\usepackage[charter]{mathdesign}
\usepackage{bm}              

\def\Rc{{\mathcal R}}
\def\Nc{{\mathcal N}}
\def\Rc{{\mathcal R}}

\def\Norm#1{\| #1 \|^2}

\def\norm#1{\| #1 \|}
\def\BT#1{\left\|#1\right\|^2}                  

\def\bF #1{\| #1 \|_F}   
\def\BF #1{\| #1 \|_F^2}  

\allowdisplaybreaks[1]  

\usepackage{graphicx}
\usepackage{subfigure}
\usepackage{epstopdf} 

\usepackage{threeparttable} 

\usepackage{enumerate}  

\usepackage{algorithmicx,algorithm, algpseudocode}

\usepackage{color}
\begin{document}

\markboth{N.-C. Wu, Q. Zuo, and T. Wang}{The deterministic single and multiple row-action methods with momentum}
\title{On the Polyak momentum variants of the greedy deterministic single and multiple row-action methods}

\author{
 Nian-Ci   Wu   \affil{1}\comma\corrauth,
 Qian      Zuo  \affil{2}\comma\affil{3}, and
 Yatian    Wang \affil{4}}

\address{
\affilnum{1}\ School of Mathematics and Statistics, South-Central Minzu University, Wuhan 430074, China.\\
\affilnum{2}\ Center on Frontiers of Computing Studies, Peking University, Beijing 10087, China.\\
\affilnum{3}\ School of Computer Science, Peking University, Beijing 10087, China.\\
\affilnum{4}\ School of Mathematics and Statistics, Wuhan University, Wuhan 430072, China.
}

\email{{\tt nianciwu@scuec.edu.cn} (N.-C. Wu)}
%
\begin{abstract}
For solving a consistent system of linear equations, the classical row-action (also known as Kaczmarz) method is a simple while really effective iteration solver. Based on the greedy index selection strategy and Polyak's heavy-ball momentum acceleration technique, we propose two deterministic row-action methods and establish the corresponding convergence theory. We show that our algorithm can linearly converge to a least-squares solution with minimum Euclidean norm. Several numerical studies have been presented to corroborate our theoretical findings.  Real-world applications, such as data fitting in computer-aided geometry design, are also presented for illustrative purposes.
\end{abstract}



\maketitle


\section{Introduction}\label{sec1}
For solving large-scale system of linear equations of the form
\begin{equation}\label{eq:Ax=b}
  Ax= b,\quad A\in \mathbb{C}^{m\times n} \quad {\rm and}\quad  b\in \mathbb{C}^{m},
\end{equation}
i.e., with $A$ being an $m$-by-$n$ complex matrix, $b$ being an $m$-dimensional complex right-hand side, and $x$ being the $n$-dimensional unknown vector, the Kaczmarz method \cite{Kac37}, also known as single row-action method, is a simple while really effective iteration solver, mainly due to its cheap per iteration cost and low total computational complexity. At each iterate of the Kaczmarz method, only one row of the coefficient matrix is used, and the new approximation is orthogonally projected onto the hyperplane which is perpendicular to this row vector. More specifically, if we use $A_{i,:}$ to represent the $i$th row of the matrix $A$ and $b_i$ the $i$th entry of the vector $b$, the Kaczmarz iteration scheme in the complex space is given as follows. For $k=0,1,2,\cdots$,
\begin{align}\label{eq:KaczIte}
   x^{(k+1)}  = x^{(k)} + \frac{b_{i_k} - A_{i_k,:} x^{(k)}}{\Norm{A_{i_k,:}}} A_{i_k,:}^\ast,
\end{align}
where the symbol $(\cdot)^\ast$ denotes the conjugate transpose of the corresponding vector or matrix, and the row index $i_k \in [m]:=\{1,2,\cdots,m\}$  is chosen according to a well-defined criterion, such as a cyclic fashion \cite{Kac37} or an appropriate probability distribution \cite{09SV}. Recently, remarkable progress of the Kaczmarz method has been made; see for example \cite{19BW, 20DSS, 212WX, 221WX, 13ZF}.

By selecting the hyperplane with a random sketch and introducing the energy norm, Gower and Richt\'{a}rik in \cite{15GR} proposed the sketch-and-project (SAP) method. In \cite{21WX}, Wu and Xiang pointed out that the SAP method is
common to the general projection methods and satisfies the Petrov-Galerkin conditions \cite{03Saad}. Gaussian Kaczmarz method is a popular representative among the SAP algorithms due to its randomness and block structure, whose iteration step in the complex space is given by
\begin{align}\label{eq:GKIte}
x^{(k+1)}  = x^{(k)} + \frac{\eta_k^\ast(b - A x^{(k)})}{\Norm{A^\ast \eta_k}} A^\ast \eta_k,
\end{align}
where $\eta_k:=[\eta_{k,i}]_{i=1}^m\in \mathbb{R}^{m}$ is a Gaussian vector with $\eta_{k,i}\sim \mathbb{N}(0,1)$. This randomized iteration scheme was analyzed in the greedy randomized average block Kaczmarz method with an adaptive stepsize \cite[Eq.(2.6)]{22MW} if we take $\eta_k = \sum_{i\in \Omega_k} (b_{i} - A_{i,:} x^{(k)}) \mu_{i}$, where $\Omega_k$ is a subset of $[m]$ and $\mu_{i}$ is the $i$th column of identity matrix with size $m$; see also \cite[Eq.(4.6)]{19Necoara}.

It is worth noting that the linear system setting described above mainly focuses on randomized iterative methods. On the other hand, it is also extremely efficient to select working rows with determinacy. McCormick used the iterative row orthogonalization to determine the non-zero vector $\eta_k$, for further details, see \cite[Eq.(9)]{77McCormick}. In a celebrated paper \cite{22CH1}, Chen and Huang updated $\eta_k$ based on a greedy criterion of row selection in \cite{18BW1,18BW2,21BW} and  provided a fast deterministic block Kaczmarz (FDBK) method.  Its relaxed version was given by \cite{22WCZ}. By applying the FDBK method on the normal equation $A^T r=0$ with $r=b-Ax$, the deterministic block coordinate descent method is obtained, see \cite{22CH2}. Shao proposed a deterministic variant of the Kaczmarz method by replacing orthogonal projection with reflection, which virtually converges more quickly than the randomized Kaczmarz method \cite{21Shao}.

Polyak momentum, popularly known as heavy ball momentum resembling the rolling of a heavy ball down the hill, is one of the oldest and influential  acceleration techniques for solving unconstrained minimization problems \cite{64Polyak} . The basic Polyak momentum update is given by
\begin{align*}
  x^{(k+1)} = x^{(k)} - \alpha \nabla f(x^{(k)}) + \beta ( x^{(k)} - x^{(k-1)} ),
\end{align*}
where $\alpha$ is a step-size, $\beta$ is a momentum parameter, and $\nabla f(x^{(k)})$ denotes the gradient of the differentiable convex function $f$ at $x^{(k)}$. When $\beta=0$, this method resolves into the so-called gradient descent method. In the context of projection-based iterative methods, the Polyak momentum technique has been incorporated into various methods, e.g., randomized coordinate descent and Kaczmarz (mRK) \cite{20LR}, sketch and project \cite{20LR}, sampling Kaczmarz Motzkin \cite{21MIA}, randomized Douglas-Rachford \cite{22HSX}, and doubly stochastic iterative framework \cite{22HX}. For other related works, we refer to \cite{12Nesterov, 16LW} and the references therein.

All existing Polyak momentum techniques for accelerating  row-action methods are randomized; see, for example \cite{20LR,21MIA,22HSX,22HX}. To the best of our knowledge, there is no convergence analysis of the deterministic and greedy row iteration schemes in the literature. It motivates us to fill this gap. In this work, blending the greedy deterministic row iteration schemes and the Polyak momentum technique, we present their Polyak momentum variant  to solve a consistent system of linear equations. The corresponding convergence theory reveals that our algorithm can converge to a least-squares solution with minimum Euclidean norm.

In the next section, we recall some basic  technical preliminaries and introduce the greedy deterministic row-action methods, including the MWRK \cite{77McCormick,19DG} and FDBK \cite{22CH1} methods. In Section \ref{Sec:Alg:mMWRK+mFDBK},  we formally present the momentum variant of the deterministic and  greedy row-action method. Their convergence theory is established in Section \ref{Sec:mMWRK+mFDBK+convergence}. In Section \ref{Sec:mMWRK+mFDBK+Numerical}, we report and discuss the numerical results. Finally, we end the paper with brief conclusions in Section \ref{Sec:mMWRK+mFDBK+Conclusions}.

{\it Notation.} Throughout the paper, all vectors are assumed to be column vectors.
For any matrix $M$,  we use $\Rc(M)$, $\Nc(M)$, $\sigma_{1}(M)$, $\sigma_{r}(M)$, $M_{i, :}$, $M_{:,j}$, and $M_{i,j}$ to denote the range, the null space, the largest,  the smallest nonzero singular values, the $i$th row, the $j$th column, and the $(i,j)$th entry,  respectively.
We remark that $M_{I_k, :}$ and $M_{:, J_k}$ stand for the row and column submatrices of $M$ indexed by index sets $I_k$ and $J_k$, respectively.
For any vector $w$,  we use $w_i$ to denote the $i$th entry.

\section{The greedy deterministic row-action methods}\label{sec:DGK}
In this section, we give a brief description of the MWRK \cite{77McCormick,19DG} and FDBK \cite{22CH1} methods, which are two typical greedy deterministic row-action methods.

A key ingredient to guarantee fast convergence of the row-action method is the construction of an appropriate criterion for the choice of row indices. In  Bai and Wu's series of works \cite{18BW1,18BW2,21BW}, the authors introduced a promising adaptive greedy index selection strategy, which is stated as follows.

{\bf The greedy row selection} \cite{18BW1,18BW2,21BW}.
{\it For $k=0,1,2,\cdots$, introducing a relaxation parameter $\theta$, the row index set is determined by
\begin{align}\label{eq:RGRS}
  U_k= \left\{i_k \in [m]\big| \psi_{i_k}(x^{(k)}) \geq
   \theta \max_{i\in [m]}\{ \psi_i(x^{(k)})\}
  + (1-\theta)\sum_{i\in [m]} \beta_i \psi_{i}(x^{(k)}) \right\}
\end{align}
with $\psi_{i}(x^{(k)}) = |b_{i} - A_{i,:}^T  x^{(k)}|^2/\Norm{A_{i,:}}$ and $\beta_i=\norm{A_{i,:}}^2 / \BF{A}$.}

Its flexibility allows us to tune the relaxation parameters and has led to several popular approaches. Specifically, the maximal weighted residual rule \cite{77McCormick, 19DG} emerges if one takes $\theta=1$, which is a generalized version of Motzkin's rule \cite{54MS} by multiplying the norm of the corresponding row of the coefficient matrix; see also \cite{21HM, 22ZL}. Later, the block version of \eqref{eq:RGRS} was given by Miao and Wu \cite{22MW}. Gower et al. enlarged this strategy in \eqref{eq:RGRS} to more general cases \cite{19GMMN}. In \cite{22HNRW}, Haddock et al. utilized a quantile of the absolute values of the residual and selected the working row. In \cite{22JZY}, Jiang et al. employed $k$-means clustering to partition the row index set. Very recently, a probability distribution depending on the angle was used to determine the row index \cite{22HDL}. We refer to \cite{21HM,17LHN,22SM,17XZ,22ZL} and the references therein for additional details on greedy row selection.

\subsection{The MWRK method}
Let $x_\ast := A^\dag b$ be a minimum Euclidean-norm least squares solution of the consistent system \eqref{eq:Ax=b}. By the Pythagorean theorem, a direct calculation indicates that the squared error of \eqref{eq:KaczIte} satisfies
\begin{align}
  \Norm{ x^{(k+1)} - x_\ast } = \Norm{ x^{(k)} - x_\ast } - \psi_{i_k}(x^{(k)}),~~i_k\in[m].
\end{align}
This implies that we may select row index such that the corresponding loss is as large as possible.  McCormick presented a deterministic greedy strategy to select the row index $i_k$ in \eqref{eq:KaczIte}, in which $i_k$ maximizes $\psi_i(x^{(k)})$ for all $i\in [m]$ \cite[Section 2.1]{77McCormick}. Accordingly, the relaxation parameter in formula \eqref{eq:RGRS} is taken by $1$. Du and Gao called it the maximal weighted residual Kaczmarz (WRRK) method, which was discussed in the real space in \cite{19DG}. Here, we reformulate it, but in the complex space, in Algorithm \ref{alg:MWRK}.

\begin{algorithm}[!htb]
\caption{The MWRK method \cite{77McCormick,19DG}.}
\label{alg:MWRK}
\begin{algorithmic}[1]
\Require
The coefficient matrix $A \in  \mathbb{C}^{m\times n}$, the right-hand side $b \in \mathbb{C}^{m}$, an initial vector $x^{(0)} \in \mathbb{C}^{n}$, and the maximum iteration number $\ell$.
\Ensure
$x^{(\ell)}$.
\State {\bf for} $k=0,1,\cdots,\ell-1$ {\bf do}
\State \quad select $i_k=\arg\max_{i\in [m]}\left\{\psi_i(x^{(k)})\right\}$;
\State \quad compute the next approximation according to \eqref{eq:KaczIte};
\State {\bf endfor}.
\end{algorithmic}
\end{algorithm}

When the initial vector is in the column space of $A^\ast$, McCormick gave an upper bound on the solution error for the MWRK iteration sequence in \cite{77McCormick}. The estimate depends on the number of equations in the system. In \cite{19DG}, Du and Gao gave a new and easily computable theoretical estimate for the convergence rate of the MWRK method. This result is restated in complex space as the following theorem.

\begin{theorem}\label{thm:MWRK} \cite[Theorem 3.1]{19DG}
Let $A\in \mathbb{C}^{m\times n}$ be a matrix without any zero rows and  $b \in \mathbb{C}^{m}$. The iteration sequence $\left\{x^{(k )}\right\}_{k=0}^{\infty}$, generated by the MWRK  method starting from any initial guess $ x^{(0)} \in \Rc(A^\ast)$, exists and  converges to the unique least-norm solution $x_{\ast}$ of the consistent linear system $Ax=b$, with the error estimate
\begin{align} \label{eq:thm-MWRK}
\Norm{x^{(k+1)} - x_{\ast}}\leq
\left(1 - \frac{\sigma_{r}^2(A)}{\widetilde{\gamma} }   \right)^k
\left(1 - \frac{\sigma_{r}^2(A)}{\BF{A}}   \right)
\Norm{x^{(0)} - x_{\ast}}
\end{align}
for $k=0,1,2,\cdots$, where the constant $\widetilde{\gamma} = \max_{i\in [m]}\left\{ \sum_{j=1,j\neq i}^{m} \norm{A_{i,:}} \right\}$.
\end{theorem}

\subsection{The FDBK method}

In this subsection, deterministic multiple row-action methods, e.g., FDBK \cite{22CH1}, are given. This approach proceeds as follows. FDBK first determines an index set using an adaptive row index selection strategy, and then applies a multiple row-action iteration scheme to update the approximation, which is listed in Algorithm \ref{alg:FDBK}.

\begin{algorithm}[!htb]
\caption{The FDBK method \cite[Algorithm 1]{22CH1}.}
\label{alg:FDBK}
\begin{algorithmic}[1]
\Require
The coefficient matrix $A \in  \mathbb{C}^{m\times n}$, the right-hand side $b \in \mathbb{C}^{m}$, an initial vector $x^{(0)} \in \mathbb{C}^{n}$, and the maximum iteration number $\ell$.
\Ensure
$x^{(\ell)}$.
\State {\bf for} $k=0,1,\cdots,\ell-1$ {\bf do}
\State \quad determine the index set $U_k$  according to formula \eqref{eq:RGRS} with $\theta=\frac{1}{2}$;
\State \quad compute $\eta_k = \sum_{i\in U_k} (b_{i} - A_{i,:} x^{(k)}) \mu_{i}$;
\State \quad compute the next approximation according to \eqref{eq:GKIte};
\State {\bf endfor}.
\end{algorithmic}
\end{algorithm}

Note that Algorithm \ref{alg:FDBK} is deterministic. It is a block Kaczmarz method, but does not require the computation of the pseudoinverses of submatrices. Chen and Huang proved that this method will converge linearly to the unique least-norm solutions of the linear systems  in the real space \cite{22CH1}. We restate this result in terms of complex space as the following.

\begin{theorem}\label{thm:FDBK} \cite[Theorem 3.1]{22CH1}
Let $A\in \mathbb{C}^{m\times n}$ be a matrix without any zero rows and  $b \in \mathbb{C}^{m}$. The iteration sequence $\left\{x^{(k )}\right\}_{k=0}^{\infty}$, generated by the FDBK  method starting from any initial guess $ x^{(0)} \in \Rc(A^\ast)$, exists and  converges to the unique least-norm solution $x_{\ast}$
of the consistent linear system $Ax=b$, with the error estimate
\begin{align} \label{eq:thm-FDBK}
\Norm{x^{(k+1)} - x_{\ast}}\leq
\left(1 - \widehat{\gamma}_k
\frac{\BF{A_{\widehat{U}_k,:}}}{\sigma_{1}^2(A_{\widehat{U}_k,:})}
\frac{\sigma_{r}^2(A)}{\BF{A}}   \right)
\Norm{x^{(k)} - x_{\ast}}
\end{align}
for $k=1,2,\cdots$,  where $\widehat{U}_k = \left\{i\in[m]| \psi_i(x^{(k)}) \neq 0 \right\}$ and
 \begin{align*}
   \widehat{\gamma}_k = \frac{1}{2}\left(
   \frac{\BF{A}}{\widehat{q}_k \BF{A_{\widehat{U}_k,:}} +(1-\widehat{q}_k) \BF{A_{U_k,:}} } +1\right)
 \end{align*}
 with $\widehat{q}_k\in(0,1]$ being a constant.
\end{theorem}

\section{The proposed methods}\label{Sec:Alg:mMWRK+mFDBK}
In this section, we propose the  momentum variant of the greedy deterministic row-action methods, including MWRK \cite{77McCormick,19DG} and FDBK \cite{22CH1}. For simplicity of notation, we name these new methods as the MWRK and FDBK methods with momentum or the mMWRK and mFDBK methods.

\subsection{The mMWRK method}

 To obtain a momentum version of the MWRK method, we use the history information to update the next iterate $x^{(k+1)}$. In particular, the mMWRK method takes two iterates $x^{(k)}$ and $x^{(k-1)}$ generated by the relaxed MWRK iteration and then updates the next iterate as a combination of the previous two updates, which is formally provided in Algorithm \ref{alg:mMWRK}.

\begin{algorithm}[!htb]
\caption{The mMWRK method.}
\label{alg:mMWRK}
\begin{algorithmic}[1]
\Require
The coefficient matrix $A \in  \mathbb{C}^{m\times n}$, the right-hand side $b \in \mathbb{C}^{m}$, two initial vectors $x^{(0)}$, $x^{(1)} \in \mathbb{C}^{n}$, a constant step-size $\alpha\in (0,2)$, a momentum parameter $\beta>0$, and the maximum iteration number $\ell$.
\Ensure
$x^{(\ell)}$.
\State {\bf for} $k=1,2,\cdots,\ell-1$ {\bf do}
\State \quad select $i_k=\arg\max_{i\in [m]}\left\{\psi_i(x^{(k)})\right\}$;
\State \quad compute the next approximation according to
\begin{align}\label{eq:mMWRKIte}
   x^{(k+1)}  = x^{(k)} + \alpha \frac{b_{i_k} - A_{i_k,:} x^{(k)}}{\Norm{A_{i_k,:}}} A_{i_k,:}^\ast + \beta (x^{(k)} - x^{(k-1)});
\end{align}
\State {\bf endfor}.
\end{algorithmic}
\end{algorithm}

\begin{figure}[!htb]
 \centering
 \subfigure[The mMWRK iteration with $(\alpha,\beta) = (1,0.5)$]{ \includegraphics[width=0.75\textwidth]{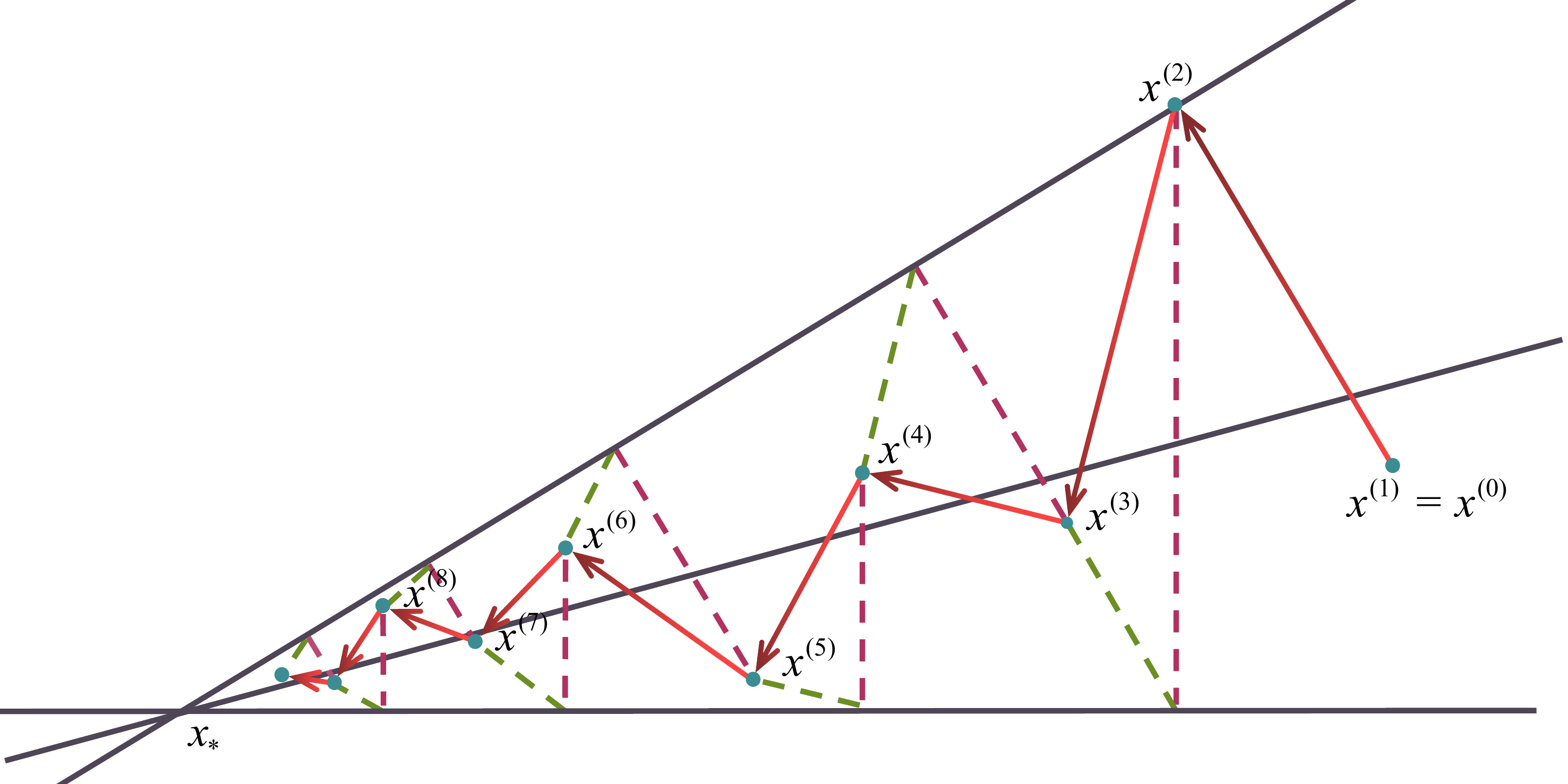}}
 \hspace{1cm}
 \subfigure[The MWRK iteration] { \includegraphics[width=0.75\textwidth]{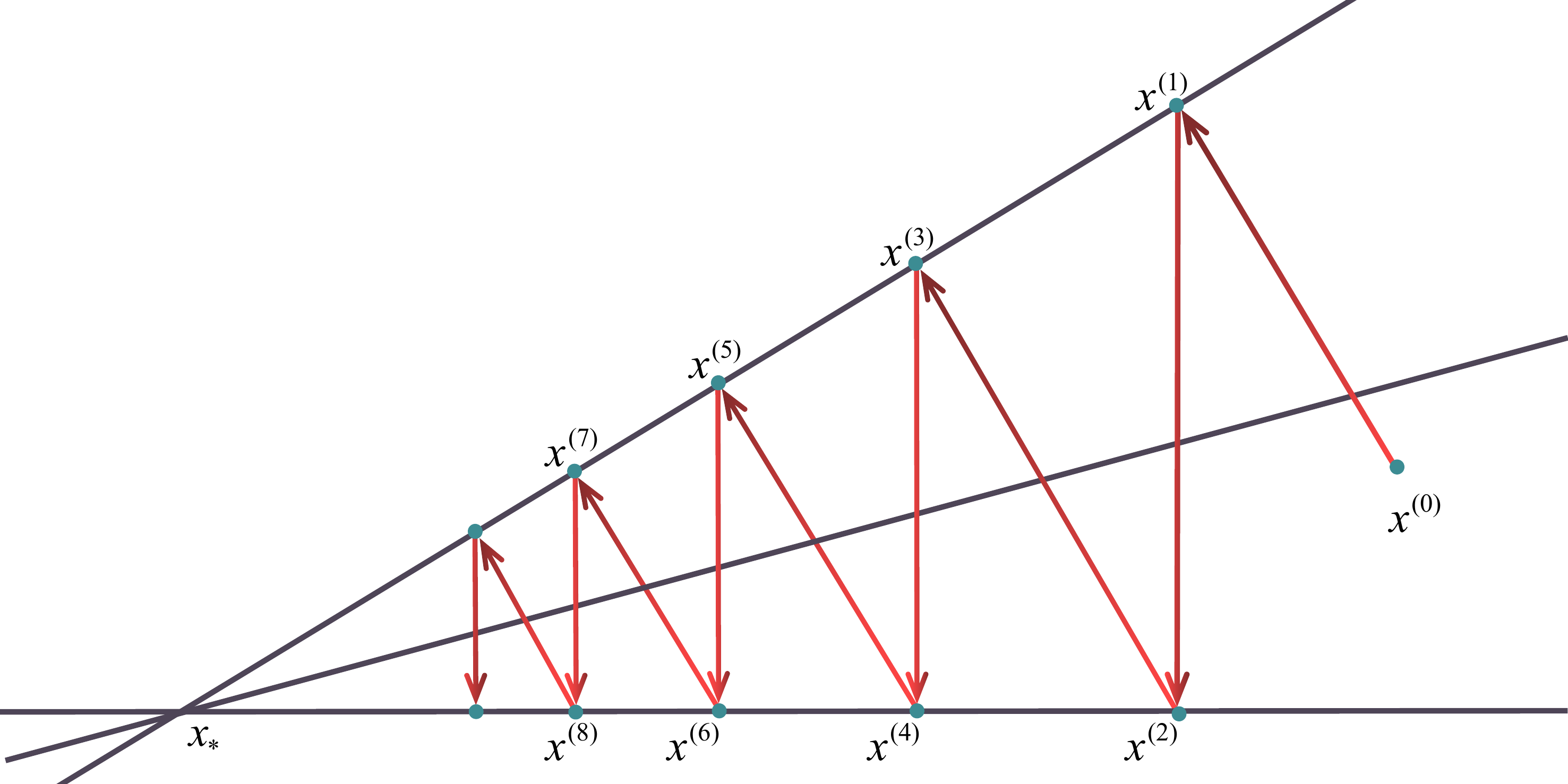}}
 \caption{ Geometric interpretation of the mMWRK (a) and MWRK (b) \cite{19DG,77McCormick} iterations. }
 \label{fig:MWRK+mMWRK-plot}
\end{figure}

\begin{remark}
To make a further insight into mMWRK, we provide a geometric interpretation of its iteration in Figure \ref{fig:MWRK+mMWRK-plot} (a). For the sake of simplicity, we consider the linear system \eqref{eq:Ax=b} with $m = 3$ and take $(\alpha,\beta) = (1,0.5)$ as an example. In such circumstances, mMWRK takes the explicit form
  \begin{align*}
    x^{(k+1)} = y^{(k)} + \frac{1}{2}(x^{(k)} - x^{(k-1)}) \quad {\rm with} \quad y^{(k)}:= x^{(k)} + \frac{b_{i_k} - A_{i_k,:} x^{(k)}}{\Norm{A_{i_k,:}}} A_{i_k,:}^\ast
  \end{align*}
for $k=1,2,\cdots$. That is, the next update $x^{(k+1)}$ (in orange-red) is obtained by forcing the momentum term $ (x^{(k)} - x^{(k-1)})/2$ (in olive-drab) to $y^{(k)}$ (in maroon), where $y^{(k)}$ is the projection of $x^{(k)}$ onto the $i_k$th hyper-plane $\{x| A_{i_k,:}x = b_{i_k}\}$ represented by a line (in gray) for $i_k = 1,2,3$ and the index $i_k$ owns the largest loss. We observe that the momentum term $(x^{(k)} - x^{(k-1)})/2$ makes the next iterate $x^{(k+1)}$ closer to the required solution ($x_{\ast}$) than that of no momentum term. Note also that the vector $x^{(k+1)}-y^{(k)}$ is always parallel to $x^{(k)} - x^{(k-1)}$ for all $k\geq 1$. With the purpose of comparison, we draw the MWRK iteration \cite{19DG,77McCormick} with the same scaling in Figure \ref{fig:MWRK+mMWRK-plot} (b). As the graph depicts,  mMWRK moves faster to $x_{\ast}$ compared to MWRK. Later, this comparison will become more apparent for larger test instances in the numerical section.
\end{remark}

\begin{remark}
The update in mMWRK is reminiscent of mRK \cite[Section 7.1]{20LR}. In fact, they are all the single row-action methods. When the index $i_k$ in \eqref{eq:mMWRKIte} is selected with  probability proportional to its Euclidean norm, the mRK method is obtained, first analyzed by Morshed et al. in \cite{20LR}. This result sparked renewed interest in the design of randomized iterative methods for solving linear systems; see, for example \cite{21MIA, 22HSX, 22HX, 12Nesterov, 16LW}.  In mMWRK, we propose to replace the selection of $i_k$ by a greedy rule. In particular, we let $i_k$ be chosen in a deterministic fashion. As far as we know, the introduction of the Polyak momentum to accelerate the deterministic and greedy single row-action methods is new.
\end{remark}

In the following, we analyze the number of flopping operations (flops) per iteration in Algorithm \ref{alg:mMWRK}. At the $k$th iterate, we first determine the index $i_k$, which needs $2m-1$ flops to compute $\psi_i(x^{(k)})$ for $i\in [m]$ and $m-1$ comparisons to obtain $\max_{i\in [m]}\left\{\psi_i(x^{(k)})\right\}$. After that, we turn to compute the next approximation.
Define the $k$th residual $r^{(k)}=b-Ax^{(k)}$ for $k=0,1,2,\cdots$. It admits the following update
\begin{align*}
 r^{(k+1)}
  = r^{(k)} - \alpha h_{k} \widetilde{A}_{:, i_k} + \beta (r^{(k)} - r^{(k-1)}),
\end{align*}
where $\widetilde{A}=AA^\ast$ and $h_{k} = r^{(k)}_{i_k}/\widetilde{A}_{i_k, i_k}$  for $i_k \in [m]$. Then, we have
\begin{align*}
 x^{(k+1)} = x^{(k)} + \alpha h_{k}~ A_{i_k,:}^\ast + \beta (x^{(k)} - x^{(k-1)}).
\end{align*}
The complexity of computing $x^{(k+1)}$ are summarized in Table \ref{alg:mMWRK-tab:ComputXk}.
\begin{table}[!htb]
\centering\renewcommand\arraystretch{1.25}
    \caption{The complexity of computing $x^{(k+1)}$ in Algorithm \ref{alg:mMWRK}.}
    \begin{tabular}{p{1.5cm}cp{7.5cm}cp{3cm}}
    \hline
\multicolumn{3}{l}{Computing~ $x^{(k+1)}$}   \\
   \hline
Step~1 && $h_k = r^{(k)}_{i_k}/\widetilde{A}_{i_k, i_k}$                                 && $1$  \\
Step~2 && $r^{(k+1)} = r^{(k)} - (\alpha h_{k}) \cdot \widetilde{A}_{:, i_k} + \beta(r^{(k)} -r^{(k-1)})$  && $5m+1$  \\
Step~3 && $x^{(k+1)} = x^{(k)} + (\alpha h_{k}) \cdot A_{i_k,:}^\ast + \beta(x^{(k)} -x^{(k-1)})$  && $5n+1$ \\
   \hline
   \end{tabular}
   \label{alg:mMWRK-tab:ComputXk}
\end{table}

\subsection{The mFDBK method}
Different from the mMWRK method choosing one index, the mFDBK method  captures all indices
in $U_k$ defined by formula \eqref{eq:RGRS}. The details of the mFDBK method are provided in Algorithm \ref{alg:mFDBK}.

\begin{algorithm}[!htb]
\caption{The mFDBK method.}
\label{alg:mFDBK}
\begin{algorithmic}[1]
\Require
The coefficient matrix $A \in  \mathbb{C}^{m\times n}$, the right-hand side $b \in \mathbb{C}^{m}$, two initial vectors $x^{(0)}$, $x^{(1)} \in \mathbb{C}^{n}$, a constant step-size $\alpha\in (0,2)$, a momentum parameter $\beta>0$, and the maximum iteration number $\ell$.
\Ensure
$x^{(\ell)}$.
\State {\bf for} $k = 1,2,\cdots,\ell-1$ {\bf do}
\State \quad determine the index set $U_k$  according to formula \eqref{eq:RGRS}  with $\theta=\frac{1}{2}$;
\State \quad compute $\eta_k = \sum_{i\in U_k} (b_{i} - A_{i,:} x^{(k)}) \mu_{i}$;
\State \quad compute the next approximation according to
 \begin{align}\label{eq:mFDBKIte}
   x^{(k+1)}  = x^{(k)} + \alpha \frac{\eta_k^\ast(b - Ax^{(k)})}{\BT{A^\ast \eta_k}} A^\ast \eta_k + \beta (x^{(k)} - x^{(k-1)});
\end{align}
\State {\bf endfor}.
\end{algorithmic}
\end{algorithm}
\begin{remark}
The acceleration mechanism behind momentum for FDBK \cite{22CH1} is similar to mMWRK for MWRK \cite{77McCormick,19DG}, the mFDBK update rule in \eqref{eq:mFDBKIte} includes two basic computational procedures. The first-half step projects $x^{(k)}$ onto the subspace $\{x| \eta_k^\ast(b - Ax)=0\}$ and the second-half step utilizes the addition of the momentum term. The value of the step-size $\alpha\in(0, 2)$ defines whether the projection is exact or not. When $\alpha \neq 1$ (with relaxation) the projection is not exact.
\end{remark}

\begin{remark}
The difference between mFDBK and FDBK \cite{22CH1} is the introduction of step-size $\alpha$ and momentum term $\beta ( x^{(k)} - x^{(k-1)})$ in the computing of $x^{(k+1)}$ in \eqref{eq:mFDBKIte}. When $(\alpha,\beta) = (1,0)$, the mFDBK method automatically reduces to the FDBK method. Whatever the parameter $\alpha$ is chosen, at each iteration step, mFDBK only requires an additional $3(m+n)$ flops.
\end{remark}

\begin{remark}\label{remark:bound+epsilon_k} The mFDBK method also inherits some basic properties of FDBK, e.g.,
\begin{enumerate}[(1)]
\setlength{\itemindent}{0cm}
\item $\BT{\eta_k}= \eta_k^\ast(b - Ax^{(k)})$;
\item $\BT{A^\ast \eta_k}\leq \sigma_{1}^2(A_{U_k,:})\BT{\eta_k}$;
\item $\epsilon_k \Norm{b - Ax^{(k)}} \leq \psi_{i} (x^{(k)})$ for any $i\in U_k$, where
$\epsilon_k =
  \max_{i\in [m]} \{ \psi_i(x^{(k)}) \}/(2\Norm{b - Ax^{(k)}}) 	
  + 1/2\BF{A}$ is a an auxiliary parameter;
\end{enumerate}
see \cite[Eqs. (3.6)-(3.9)]{22CH1}. In addition, from the definition of $\epsilon_k$, we can draw a conclusion that
\begin{align}\label{eq:bound_epsilon}
\epsilon_k \BF{A}
& = \frac{1}{2}\frac{\BF{A}}{\Norm{b - Ax^{(k)}}}
  \max_{i\in [m]}\left\{ \psi_i(x^{(k)}) \right\}
  + \frac{1}{2}\notag\\
& = \frac{1}{2}\frac{\BF{A}}
{ \sum_{i\in\widehat{U}_k} \psi_i(x^{(k)})\Norm{A_{i,:}}}
  \max_{i\in [m]}\left\{ \psi_i(x^{(k)}) \right\}
  + \frac{1}{2}\notag\\
&\geq \frac{1}{2} \frac{\BF{A}}
{ \BF{A_{\widehat{U}_k,:}}}
  + \frac{1}{2}.
\end{align}
\end{remark}

The mFDBK method is an adaptive deterministic multiple row-action method. The index set $U_k$ is updated at each step as the iteration proceeds. It can be implemented at low computational cost and analyzed in a relatively simple way as follows.

At step $4$ in Algorithm \ref{alg:mFDBK}, we have
\begin{align*}
r^{(k+1)}
 & =r^{(k)} -
 \alpha\frac{\sum_{i\in U_k} r^{(k)}_{i} \mu_{i}^Tr^{(k)}}
      {\sum_{i\in U_k} \sum_{j\in U_k} r^{(k)}_{i}  r^{(k)}_{j} \mu_{i} AA^\ast \mu_{j}^T}
       \sum_{i\in U_k}  r^{(k)}_{i} AA^\ast \mu_{i} + \beta(r^{(k)} -r^{(k-1)})\\
 & =r^{(k)} -
 \alpha\frac{\sum_{i\in U_k}\left( r^{(k)}_{i}\right)^2}
      {\sum_{i\in U_k} \sum_{j\in U_k}  r^{(k)}_{i} r^{(k)}_{j} \widetilde{A}_{i, j} }
       \sum_{i\in U_k}  r^{(k)}_{i} \widetilde{A}_{:,i}+ \beta(r^{(k)} -r^{(k-1)}).
\end{align*}
Then, the  next approximation is computed by
\begin{align*}
 x^{(k+1)}
  = x^{(k)} +
 \alpha\frac{\sum_{i\in U_k}\left( r^{(k)}_{i}\right)^2}
      {\sum_{i\in U_k} \sum_{j\in U_k}  r^{(k)}_{i} r^{(k)}_{j} \widetilde{A}_{i, j} }
       \sum_{i\in U_k}  r^{(k)}_{i} A_{i,:}^\ast + \beta(x^{(k)} -x^{(k-1)}).
\end{align*}
Based on the recursive update formula of $x^{(k)}$, Algorithm \ref{alg:mFDBK} really only updates two vectors $r^{(k)}$ and $x^{(k)}$. We emphasize that the blocking strategy is always applied to the matrix-vector multiplication, e.g.,
\begin{align*}
   \sum_{i\in U_k}  r^{(k)}_{i} \widetilde{A}_{:,i}
   = \widehat{A}_{:,U_k}r^{(k)}_{U_k}
  ~~{\rm and}~~
   \sum_{i\in U_k} \sum_{j\in U_k}  r^{(k)}_{i} r^{(k)}_{j} \widetilde{A}_{i, j}
   = \left( r^{(k)}_{U_k}\right)^T  \left( \widetilde{A}_{U_k, U_k}   r^{(k)}_{U_k}    \right).
\end{align*}
The computational procedure is summarized in Table \ref{alg:mFDBK-tab:ComputXk}. Note that part of the component at iteration $k$ of mFDBK is the construction of the index set $U_k$. We recognize that if we have $r^{(k)}$ and $x^{(k)}$ at the beginning, the resulting method could be made faster. For more details of computing $U_k$, we refer to \cite{18BW1,18BW2,21BW}.

\begin{table}[!htb]
\centering\renewcommand\arraystretch{1.25}
    \caption{The complexity of computing $x^{(k+1)}$ in Algorithm \ref{alg:mFDBK}.}
    \begin{tabular}{p{1.5cm}cp{7.5cm}cp{3cm}}
    \hline
\multicolumn{3}{l}{Computing~ $x^{(k+1)}$}   \\
   \hline
Step~1 && $\widetilde{g}_{k} = \left(r^{(k)}_{U_k}\right)^\ast \left(r^{(k)}_{U_k}\right)$  && $2|U_k|-1$ \\
Step~2 && $\widehat{g}_{k} = \left(r^{(k)}_{U_k}\right)^\ast
   \left( \widetilde{A}_{U_k, U_k}r^{(k)}_{U_k} \right)$  &&$ 3(|U_k|^2+|U_k|)/2 $  \\
Step~3 && $g_{k} = \widetilde{g}_{k}/\widehat{g}_{k}$                                 && $1$  \\
Step~4 && $g^{(k,1)} = \widetilde{A}_{:,U_k}r^{(k)}_{U_k}$   && $m(2|U_k|-1)$  \\
Step~5 && $r^{(k+1)} = r^{(k)} - (\alpha ~g_{k}) \cdot g^{(k,1)} + \beta(r^{(k)} -r^{(k-1)})$  && $5m+1$  \\
Step~6 && $g^{(k,2)} = A_{U_k,:}^\ast r^{(k)}_{U_k}$  && $n(2|U_k|-1)$  \\
Step~7 && $x^{(k+1)} = x^{(k)} + (\alpha ~g_{k}) \cdot g^{(k,2)} + \beta(x^{(k)} -x^{(k-1)})$  && $5n+1$ \\
   \hline
   \end{tabular}
   \begin{tablenotes}
   \footnotesize
   \item[1.] {\it The symbol $|U|$ denotes the cardinality of a set $U$.}
   \end{tablenotes}
   \label{alg:mFDBK-tab:ComputXk}
\end{table}

\section{Convergence analyses}\label{Sec:mMWRK+mFDBK+convergence}

The aim of this section is to demonstrate the convergence properties of the mMWRK and mFDBK methods for solving a consistent system of linear equations. First we present a lemma from \cite{20LR} which we will use in our convergence proof.

\begin{lemma}\label{mFDBK:Fk_123} \cite[Lemma 9]{20LR}
  Fix $F_1 = F_0 \geq 0$, $\zeta \geq 0$, and let $\left\{F_k\right\}_{k=0}^{\infty}$ be a sequence of nonnegative real
numbers satisfying the relation
\begin{align*}
  F_{k+1} \leq a_1 F_{k} + a_2 F_{k-1}
\end{align*}
for any $k \geq 1$, where $a_2 \geq 0$, $a_1 + a_2 <1$. Then the sequence satisfies the relation
\begin{align*}
   F_{k+1} \leq q^k (1+p)F_{0}
\end{align*}
for any $k \geq 0$, where $p = (\sqrt{a_1^2 + 4 a_2} - a_1)/2$ and $q = (\sqrt{a_1^2 + 4 a_2} + a_1)/2$.
\end{lemma}

A proof of this lemma can be found in \cite{20LR}.

\begin{theorem}\label{thm:mMWRK}
 Let the linear system \eqref{eq:Ax=b}, with the coefficient matrix $A\in \mathbb{C}^{m\times n}$ and the right-hand side $b \in \mathbb{C}^{m}$, be consistent.
 The mMWRK method (see Algorithm \ref{alg:mMWRK}) starts from any initial guesses $ x^{(0)} = x^{(1)}$ in the column space of $A^\ast$ and generates the iteration sequence $\left\{x^{(k)}\right\}_{k=0}^{\infty}$, then the next squared error satisfies
\begin{align*}
  \Norm{x^{(k+1)}-x_\ast } \leq \gamma_{k,1} \Norm{ x^{(k)}-x_\ast } + \gamma_{k,2} \Norm{ x^{(k-1)}-x_\ast }
\end{align*}
for $k=1,2,\cdots$, where
\begin{align*}
 \gamma_{k,1} = (1+3\beta+\beta^2) + (\alpha^2 - 2\alpha - \alpha\beta)\widetilde{\rho}_k
 \quad {\rm and} \quad
 \gamma_{k,2} = 2\beta^2+(1+\alpha)\beta
\end{align*}
with $\widetilde{\rho}_k = \sigma_r^2(A)/ \Norm{A_{\widehat{U}_k},:}_F$ and $\widehat{U}_k = \left\{i\in[m]| \psi_i(x^{(k)}) \neq 0 \right\}$.
\end{theorem}

\begin{proof}
At the $(k+1)$th mMWRK iteration,  by substituting the update \eqref{eq:mMWRKIte} into the next squared error for $k=0,1,2,\cdots$, we divide it into three parts, i.e.,
\begin{align}\label{eq1:mMWRK_x_k+1}
  \Norm{ x^{(k+1)}-x_\ast } = s_{k,1}  + s_{k,2} + s_{k,3},
\end{align}
where $s_{k,1}$,  $s_{k,2}$,  and  $s_{k,3}$ are respectively defined by
\begin{align*}
\left \{
\begin{array}{l}
   s_{k,1}  =  \Norm{x^{(k)} - x_{\ast} + \alpha  h_{k}  A_{i_k,:}^{\ast}},  \vspace{1.5ex}\\
   s_{k,2}  =  2\beta\left\langle x^{(k)} - x_{\ast} + \alpha h_{k} A_{i_k,:}^{\ast},~~ x^{(k)} - x^{(k-1)} \right\rangle,  \vspace{1.5ex} \\
   s_{k,3} = \beta^2 \Norm{ x^{(k)} - x^{(k-1)} }.
\end{array}
\right.
\end{align*}
We proceed to analyze them individually.

According to the fact that $\Norm{h_{k} A_{i_k,:}} = \psi_{i_k}(x^{(k)})$, we have
\begin{align}\label{eq:mMWRK_sk1}
  s_{k,1}
  & = \Norm{x^{(k)} - x_{\ast}}
      + 2\alpha \left\langle x^{(k)} - x_{\ast},~~ h_{k}  A_{i_k,:}^{\ast} \right\rangle
      + \alpha^2 \Norm{h_{k} A_{i_k,:}}\notag\\
  & = \Norm{x^{(k)} - x_{\ast}} + (\alpha^2 - 2\alpha) \psi_{i_k}(x^{(k)}).
\end{align}
Define two auxiliary variables
\begin{align*}
  s_{k,2}^{(1)}
  & := 2\beta\left\langle x^{(k)} - x_{\ast},~~ x^{(k)}- x_{\ast} \right\rangle  +
      2\beta\left\langle x^{(k)} - x_{\ast},~~ x_{\ast}  - x^{(k-1)} \right\rangle \\
  &\leq 2\beta\Norm{x^{(k)} - x_{\ast}} +  \beta (\Norm{x^{(k)} - x_{\ast}} + \Norm{x_{\ast}  - x^{(k-1)}})\\
  & = 3\beta\Norm{x^{(k)} - x_{\ast}} + \beta \Norm{ x^{(k-1)} - x_{\ast}}\\
  s_{k,2}^{(2)}
  & := 2\alpha\beta\left\langle h_{k} A_{i_k,:}^{\ast},~~ x^{(k)}- x_{\ast} \right\rangle +
      2\alpha\beta\left\langle h_{k} A_{i_k,:}^{\ast},~~ x_{\ast}  - x^{(k-1)} \right\rangle\\
  & =-2\alpha\beta \psi_{i_k}(x^{(k)})   +
      2\alpha\beta\left\langle h_{k} A_{i_k,:}^{\ast},~~ x_{\ast}  - x^{(k-1)} \right\rangle\\
  &\leq-2\alpha\beta \psi_{i_k}(x^{(k)})   +
       \alpha\beta (\psi_{i_k}(x^{(k)}) + \Norm{x^{(k-1)} - x_{\ast}})\\
  & =  \alpha\beta \Norm{x^{(k-1)} - x_{\ast}} -  \alpha\beta \psi_{i_k}(x^{(k)}).
\end{align*}
It follows that
\begin{align}\label{eq:mMWRK_sk2}
  s_{k,2}
& = 2\beta\left\langle x^{(k)} - x_{\ast},~~ x^{(k)} - x^{(k-1)} \right\rangle +
    2\alpha\beta\left\langle h_{k} A_{i_k,:}^{\ast},~~ x^{(k)} - x^{(k-1)} \right\rangle\notag\\
&:=s_{k,2}^{(1)}  +   s_{k,2}^{(2)}\notag\\
& \leq 3\beta\Norm{x^{(k)} - x_{\ast}} + (1+\alpha)\beta \Norm{ x^{(k-1)} - x_{\ast}} -  \alpha\beta \psi_{i_k}(x^{(k)}).
\end{align}
Using the inequality $\Norm{x-y}\leq 2(\Norm{x-z} + \Norm{y-z})$ for any vectors $x$, $y$, and $z$ with compatible dimension, it holds that
\begin{align}\label{eq:mMWRK_sk3}
 s_{k,3} \leq 2\beta^2 \Norm{ x^{(k)} - x_\ast } + 2\beta^2 \Norm{x^{(k-1)}- x_\ast }.
\end{align}

Combining formulas \eqref{eq:mMWRK_sk1}, \eqref{eq:mMWRK_sk2}, and \eqref{eq:mMWRK_sk3}, it indicates that
\begin{align*}
  \Norm{ x^{(k+1)} - x_\ast }
  &\leq (1+3\beta+\beta^2)\Norm{ x^{(k)} - x_\ast }
  + (\alpha^2 - 2\alpha - \alpha\beta) \psi_{i_k}(x^{(k)})\\
  &\quad + (2\beta^2+(1+\alpha)\beta)\Norm{x^{(k-1)}- x_\ast }.
\end{align*}
Since $\psi_{i_k}(x^{(k)}) = \max_{i\in[m]}\left\{\psi_i(x^{(k)})\right\}$ at the $k$th iterate, it follows that
\begin{align*}
  \psi_{i_k}(x^{(k)})
& = \psi_{i_k}(x^{(k)}) \frac{\Norm{b-Ax^{(k)}}}{ \sum_{i=1}^m \psi_i(x^{(k)}) \Norm{A_{i,:}} } \\
& = \psi_{i_k}(x^{(k)}) \frac{\Norm{b-Ax^{(k)}}}{ \sum_{i\in \widehat{U}_k} \psi_i(x^{(k)}) \Norm{A_{i,:}} } \\
& \geq \frac{\Norm{b-Ax^{(k)}}}{ \Norm{A_{\widehat{U}_k},:}_F } \\
&  \geq \widetilde{\rho}_k \Norm{ x^{(k)}-x_\ast }.
\end{align*}
Here in the last inequality we have used the estimate
\begin{align}\label{eq:Au=sigma}
  \Norm{A (x^{(k)}-x_\ast) }\geq \sigma_r^2(A)\Norm{ x^{(k)}-x_\ast},
\end{align}
which holds true because $ x^{(k)}-x_\ast $ belongs to the column space of $A^\ast$ shown by induction. Then, we straightforwardly obtain the estimate \eqref{eq1:mMWRK_x_k+1}.
\qed
\end{proof}

By utilizing the fact that
\begin{align*}
  1\geq\widetilde{\rho}_k = \frac{\Norm{A}_F}{\Norm{A_{\widehat{U}_k},:}_F}  \frac{\sigma_r^2(A)}{\Norm{A}_F} > \rho>0,
\end{align*}
the mMWRK squared error satisfies that
\begin{align*}
  \Norm{ x^{(k+1)}-x_\ast } \leq \gamma_{1} \Norm{ x^{(k)}-x_\ast } + \gamma_{2} \Norm{ x^{(k-1)} -x_\ast }
\end{align*}
with $\gamma_{1}=(1+3\beta+\beta^2) + (\alpha^2 - 2\alpha - \alpha\beta)\rho$ and $\gamma_{2} = \gamma_{k,2}>0$ for $k=1,2,\cdots$.
Then, we apply Lemma \ref{mFDBK:Fk_123}, wherein the two coefficients are given above, and obtain a convergence result stated as follows.

\begin{remark}\label{Remark:gamma1+gamma2<1}
Let $f(\beta)=3\beta^2+\tau_1 \beta+\tau_2$, where $\tau_1=4+\alpha-\alpha \rho$ and $\tau_2=\alpha(2-\alpha)\rho$. The estimate $f(\beta)<0$ is true for
\begin{align*}
  0<\alpha<2 \quad and \quad
  0<\beta<(\sqrt{\tau_1^2 + 12\tau_2}-\tau_1)/6.
\end{align*}
It implies that
\begin{align*}
  \gamma_{1}+\gamma_{2} = f(\beta)+1<1.
\end{align*}
Thus, the assumption for Lemma \ref{mFDBK:Fk_123} holds, so we have that
\begin{align*}
  \Norm{ x^{(k+1)} - x_\ast} \leq
  \left( \frac{\sqrt{\gamma_1^2 + 4 \gamma_2} + \gamma_1}{2}\right)^k
  \left( \frac{\sqrt{\gamma_1^2 + 4 \gamma_2} - \gamma_1}{2}\right)
  \Norm{ x^{(0)}  - x_\ast }.
\end{align*}
The convergence factor is less than $1$ following directly from the assumption $\gamma_{1}+\gamma_{2} < 1$ since
\begin{align*}
\frac{\sqrt{\gamma_1^2 + 4 \gamma_2} + \gamma_1}{2} - 1 =
\frac{\sqrt{(\gamma_1-2)^2 + 4 (\gamma_1+\gamma_2)-4} + (\gamma_1 - 2)}{2}<0.
\end{align*}
\end{remark}

Now, we turn to analyze the convergence of the mFDBK method.

\begin{theorem}\label{thm:mFDBK}
Let the linear system \eqref{eq:Ax=b}, with the coefficient matrix $A\in \mathbb{C}^{m\times n}$ and the right-hand side $b \in \mathbb{C}^{m}$, be consistent.
 The mFDBK method (see Algorithm \ref{alg:mFDBK}) starts from any initial guesses $ x^{(0)} = x^{(1)}$ in the column space of $A^\ast$ and generates the iteration sequence $\left\{x^{(k )}\right\}_{k=0}^{\infty}$, then
the next squared error satisfies
\begin{align*}
  \Norm{ x^{(k+1)}-x_\ast } \leq \gamma_{k,3} \Norm{ x^{(k)}-x_\ast } + \gamma_{k,2} \Norm{ x^{(k-1)}-x_\ast }
\end{align*}
for $k=1,2,\cdots$, where
\begin{align*}
 \gamma_{k,3} = (1+3\beta+\beta^2) + (\alpha^2 - 2\alpha - \alpha\beta)\widehat{\rho}_k
\end{align*}
with
\begin{align*}
\widehat{\rho}_k = \left( \frac{1}{2} \frac{\BF{A}}
{ \BF{A_{\widehat{U}_k,:}}}
  + \frac{1}{2} \right)
  \frac{ \Norm{A_{U_k,:}}_F }{\BF{A}}
\frac{\sigma_{r}^2(A)}{ \sigma_{1}^2(A_{U_k,:})}
\end{align*}
and $\gamma_{k,2}$ and $\widehat{U}_k$ being defined by Theorem \ref{thm:mMWRK}.
\end{theorem}

\begin{proof}
The proof is similar to that of Theorem \ref{thm:mMWRK} in,  with slightly different technicalities involved. For completeness and simplicity, we first write the $(k+1)$th mFDBK squared error as follows.
\begin{align}\label{eq1:mFDBK_x_k+1}
  \Norm{ x^{(k+1)} -x_\ast } = s_{k,4}  + s_{k,5} + s_{k,6},
\end{align}
where $s_{k,4}$,  $s_{k,5}$,  and  $s_{k,6}$ are respectively defined by
\begin{align*}
\left \{
\begin{array}{l}
   s_{k,4}  =  \Norm{x^{(k)} - x_{\ast} + \alpha  g_{k}  A^\ast \eta_k},  \vspace{1.5ex}\\
   s_{k,5}  =  2\beta\left\langle x^{(k)} - x_{\ast} + \alpha g_{k} A^\ast \eta_k,~~ x^{(k)} - x^{(k-1)} \right\rangle,  \vspace{1.5ex} \\
   s_{k,6} = \beta^2 \Norm{ x^{(k)} - x^{(k-1)} }.
\end{array}
\right.
\end{align*}
These three terms will be analyzed individually.

According to the fact that
\begin{align*}
  \left\langle x_{\ast} - x^{(k)},~~ g_{k}  A^\ast \eta_k \right\rangle
  = \Norm{g_{k} A^\ast \eta_k}
  = \phi(x^{(k)}, \eta_k)
\end{align*}
with $\phi(x^{(k)}, \eta_k)=  | \eta_k^T(b - Ax^{(k)}) |^2/\Norm{A^T \eta_k}$, we rewrite the first term in \eqref{eq1:mFDBK_x_k+1} as
\begin{align}\label{eq:mFDBK_tk1}
  s_{k,4}
  & = \Norm{x^{(k)} - x_{\ast}}
      + 2\alpha \left\langle x^{(k)} - x_{\ast},~~ g_{k}  A^\ast \eta_k \right\rangle
      + \alpha^2 \Norm{g_{k} A^\ast \eta_k}\notag\\
  & = \Norm{x^{(k)} - x_{\ast}} + (\alpha^2 - 2\alpha) \phi(x^{(k)}, \eta_k).
\end{align}
Define two auxiliary variables
\begin{align*}
  s_{k,5}^{(1)}
  & := 2\beta\left\langle x^{(k)} - x_{\ast},~~ x^{(k)}- x_{\ast} \right\rangle  +
      2\beta\left\langle x^{(k)} - x_{\ast},~~ x_{\ast}  - x^{(k-1)} \right\rangle \\
  &\leq 2\beta\Norm{x^{(k)} - x_{\ast}} +  \beta (\Norm{x^{(k)} - x_{\ast}} + \Norm{x_{\ast}  - x^{(k-1)}})\\
  & = 3\beta\Norm{x^{(k)} - x_{\ast}} + \beta \Norm{ x^{(k-1)} - x_{\ast}}\\
  s_{k,5}^{(2)}
  & := 2\alpha\beta\left\langle g_{k} A^\ast \eta_k,~~ x^{(k)}- x_{\ast} \right\rangle +
      2\alpha\beta\left\langle g_{k} A^\ast \eta_k,~~ x_{\ast}  - x^{(k-1)} \right\rangle\\
  & =-2\alpha\beta \phi(x^{(k)}, \eta_k)   +
      2\alpha\beta\left\langle g_{k} A^\ast \eta_k,~~ x_{\ast}  - x^{(k-1)} \right\rangle\\
  &\leq-2\alpha\beta \phi(x^{(k)}, \eta_k)   +
       \alpha\beta \phi(x^{(k)}, \eta_k) + \alpha\beta \Norm{x^{(k-1)} - x_{\ast}} \\
  & =  \alpha\beta \Norm{x^{(k-1)} - x_{\ast}} -  \alpha\beta \phi(x^{(k)}, \eta_k).
\end{align*}
It follows that
\begin{align}\label{eq:mFDBK_tk2}
  s_{k,5}
& = 2\beta\left\langle x^{(k)} - x_{\ast},~~ x^{(k)} - x^{(k-1)} \right\rangle +
    2\alpha\beta\left\langle g_{k} A^\ast \eta_k,~~ x^{(k)} - x^{(k-1)} \right\rangle\notag\\
&:=s_{k,5}^{(1)}  +   s_{k,5}^{(2)}\notag\\
& \leq 3\beta\Norm{x^{(k)} - x_{\ast}} + (1+\alpha)\beta \Norm{ x^{(k-1)} - x_{\ast}} -  \alpha\beta \phi(x^{(k)}, \eta_k).
\end{align}
By adding and subtracting $x_\ast$ for the third term of formula \eqref{eq1:mFDBK_x_k+1},
\begin{align*}
s_{k,6} = \beta^2 \Norm{ (x^{(k)} - x_\ast ) + (x_\ast - x^{(k-1)}) },
\end{align*}
we have
\begin{align}\label{eq:mFDBK_tk3}
 s_{k,6} \leq 2\beta^2 \Norm{ x^{(k)} - x_\ast } + 2\beta^2 \Norm{x^{(k-1)}- x_\ast }.
\end{align}

Grouping like terms in formulas \eqref{eq:mFDBK_tk1}, \eqref{eq:mFDBK_tk2}, and \eqref{eq:mFDBK_tk3}, it indicates that
\begin{align*}
  \Norm{ x^{(k+1)} - x_\ast }
  &\leq (1+3\beta+\beta^2)\Norm{ x^{(k)} - x_\ast }
  + (\alpha^2 - 2\alpha - \alpha\beta) \phi(x^{(k)}, \eta_k)\\
  &\quad + (2\beta^2+(1+\alpha)\beta)\Norm{x^{(k-1)}- x_\ast }.
\end{align*}
Therefore, the conclusion follows by utilizing the fact that
\begin{align*}
\phi(x^{(k)}, \eta_k)
& = \frac{\sum_{i_k\in U_k} \big| b_{i_k} - A_{i_k,:}^T x^{(k)} \big|^2 \cdot \Norm{\eta_k}}{\Norm{A^T \eta_k}}\\
&\geq  \frac{\sum_{i_k\in U_k} \psi_{i_k}(x^{(k)}) \Norm{A_{i_k,:}}}{ \sigma_{1}^2(A_{U_k,:})}\\
&\geq \frac{\sum_{i_k\in U_k} \epsilon_k \Norm{b - Ax^{(k)}}\Norm{A_{i_k,:}}}{ \sigma_{1}^2(A_{U_k,:})}\\
&= \frac{\epsilon_k \Norm{b - Ax^{(k)}} \Norm{A_{U_k,:}}_F }{ \sigma_{1}^2(A_{U_k,:})}\\
&\geq \epsilon_k \BF{A}
\frac{ \Norm{A_{U_k,:}}_F }{\BF{A}}
\frac{\sigma_{r}^2(A)}{ \sigma_{1}^2(A_{U_k,:})}
\Norm{x^{(k)} - x_\ast }\\
& \geq\left( \frac{1}{2} \frac{\BF{A}}
{ \BF{A_{\widehat{U}_k,:}}}
  + \frac{1}{2} \right)
  \frac{ \Norm{A_{U_k,:}}_F }{\BF{A}}
\frac{\sigma_{r}^2(A)}{ \sigma_{1}^2(A_{U_k,:})}
\Norm{x^{(k)}  - x_\ast },
\end{align*}
where the first and second inequalities are from Remark \ref{remark:bound+epsilon_k}, the third inequality is from formula \eqref{eq:Au=sigma} since $x^{(k)} - x_\ast$ is in the column space of $A^\ast$ by induction, and the last inequality follows from formula \eqref{eq:bound_epsilon}.
\qed
\end{proof}

\begin{remark}
We remain to check that $0 < \rho \leq \widehat{\rho}_k < 1$ for $k = 1,2,\cdots$. It is a similar story to obtain the convergence result in Remark
 \ref{Remark:gamma1+gamma2<1}. The details here are omitted.
\end{remark}

\section{Numerical simulations}\label{Sec:mMWRK+mFDBK+Numerical}

In this section, we implement the Polyak momentum variant of the greedy deterministic row-action methods (e.g., mMWRK and mFDBK) and its original variant (e.g., MWRK \cite{77McCormick,19DG} and FDBK \cite{22CH1}),  and show that the former is numerically advantageous over the latter in terms of the number of iteration steps, which is abbreviated as IT. We also report the speed-up (SU) of mMWRK against MWRK and mFDBK against FDBK, which are defined as
\begin{align*}
  {\rm SU_1} = \frac{{\rm IT~of~MWRK}}{{\rm IT~of~mMWRK}}
  \quad {\rm and}\quad
  {\rm SU_2} = \frac{{\rm IT~of~FDBK}}{{\rm IT~of~mFDBK}}.
\end{align*}

In our implementations, one of the solution vectors is set by $x_\ast = A^\dag e \in \mathbb{R}^n$, where the vector $e$ is generated by using the MATLAB function {\tt ones}$(m,1)$, and the right-hand side $b \in \mathbb{C}^m$ is taken to be $Ax_\ast$. All computations are respectively started from $x^{(0)}=x^{(1)}=\textsc{0}$ and $x^{(0)} =\textsc{0}$ with and without momentum accelerations, and terminated once the relative solution error (RSE), defined by RSE $= \Norm{x^{(k)} - x_{\ast}}/\Norm{x_{\ast}}$ at the current iterate $x^{(k)}$, satisfies ${\rm RSE}\leq 10^{-12}$, or the number of iteration steps exceeds $10^5$. In addition, we execute mMWRK, mFDBK, MWRK, and FDBK without explicitly forming the matrices $\widetilde{A}=AA^\ast$.  All numerical tests are performed on a Founder desktop PC with Intel(R) Core(TM) i5-7500 CPU 3.40 GHz.

We note that in the single row-action methods, the authors in \cite{19DG} have shown that the MWRK method is more efficient than the greedy randomized Kaczmarz in \cite{18BW1}. Numerical results in \cite{22CH1}  illustrated that the FDBK method provides more significant computational advantages than several existing multiple row-action methods, including the randomized average block Kaczmarz method \cite{19Necoara} and the greedy block Kaczmarz method \cite{20NZ}. Then, in the following, we just compare the efficiency of mMWRK (resp. mFDBK) with MWRK (resp. FDBK).

\subsection{Choice of $\alpha$ and $\beta$}\label{sec:choice_alpha+beta}
In this subsection, we demonstrate the computational behavior of mMWRK and mFDBK with respect to different step-sizes $\alpha$ and momentum parameters $\beta$. The coefficient matrix is the synthetic data, generated by the MATLAB function {\tt randn}$(m,n)$.

\begin{figure}[!htb]
\centering
    \subfigure[mMWRK: $n=50$ and $m=10n$]{
		\includegraphics[width=0.48\textwidth]{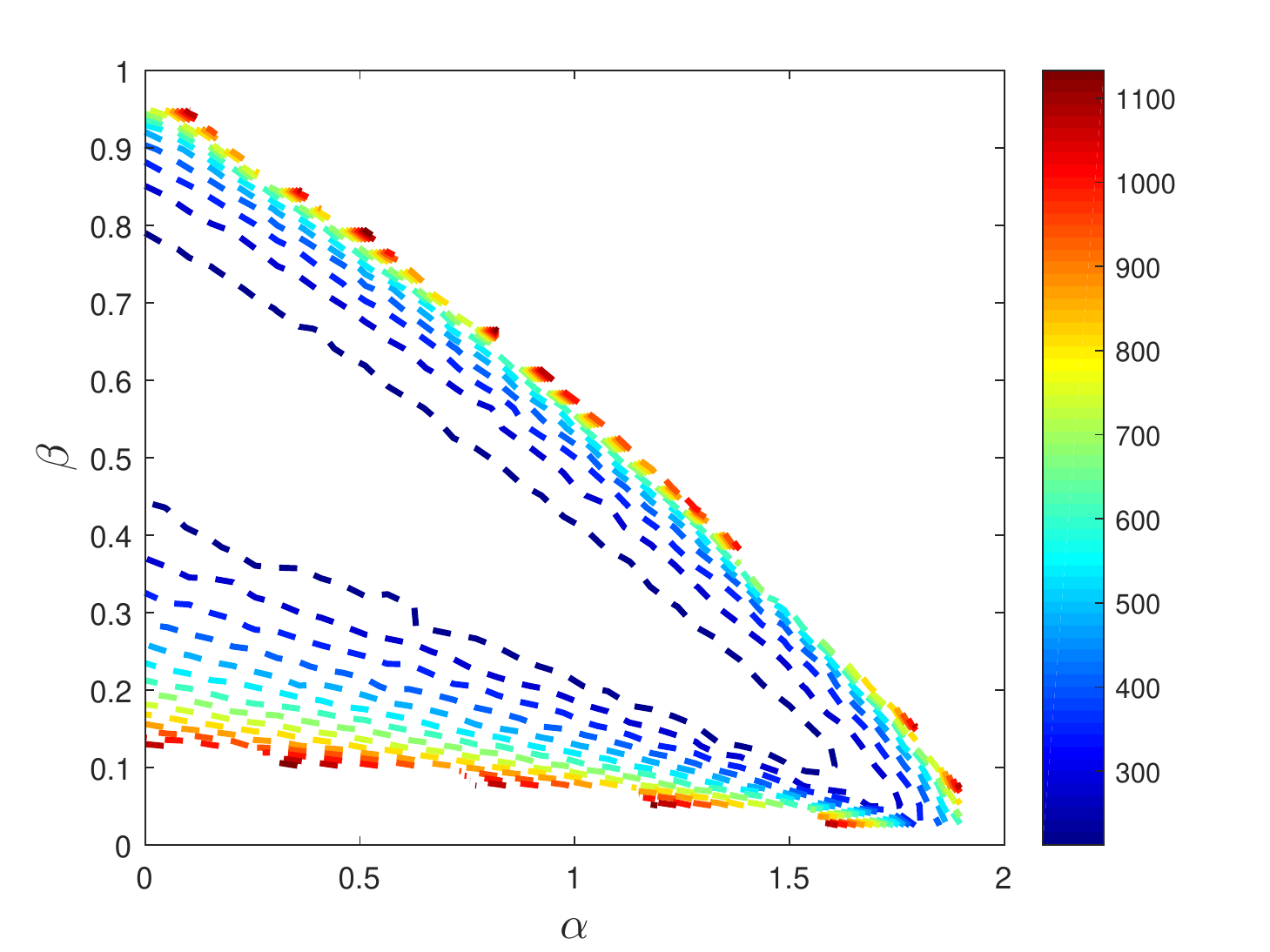}}
	\subfigure[mFDBK: $n=50$ and $m=10n$]{
	    \includegraphics[width=0.48\textwidth]{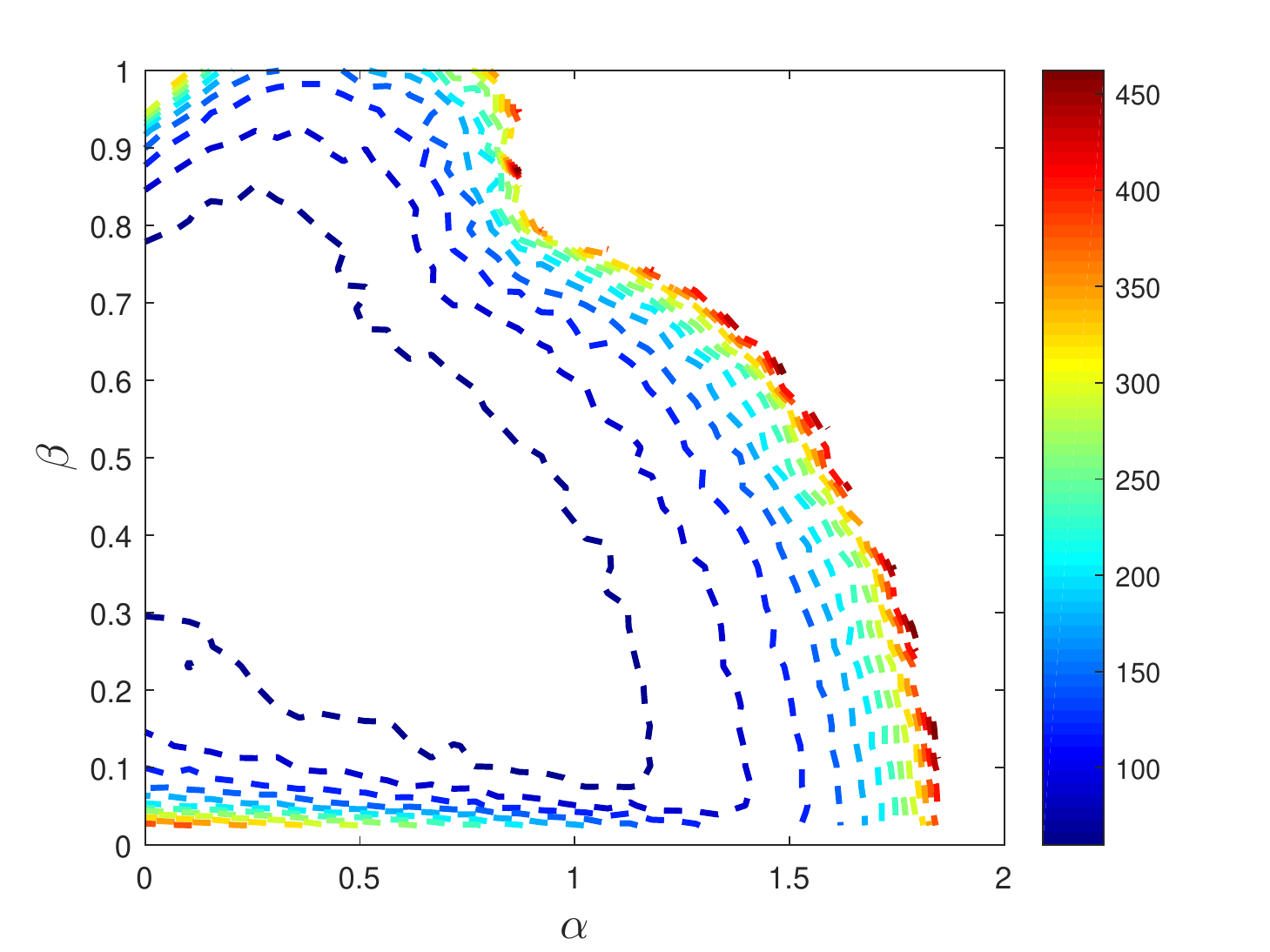}}
    \subfigure[mMWRK: $n=50$ and $m=50n$]{
		\includegraphics[width=0.48\textwidth]{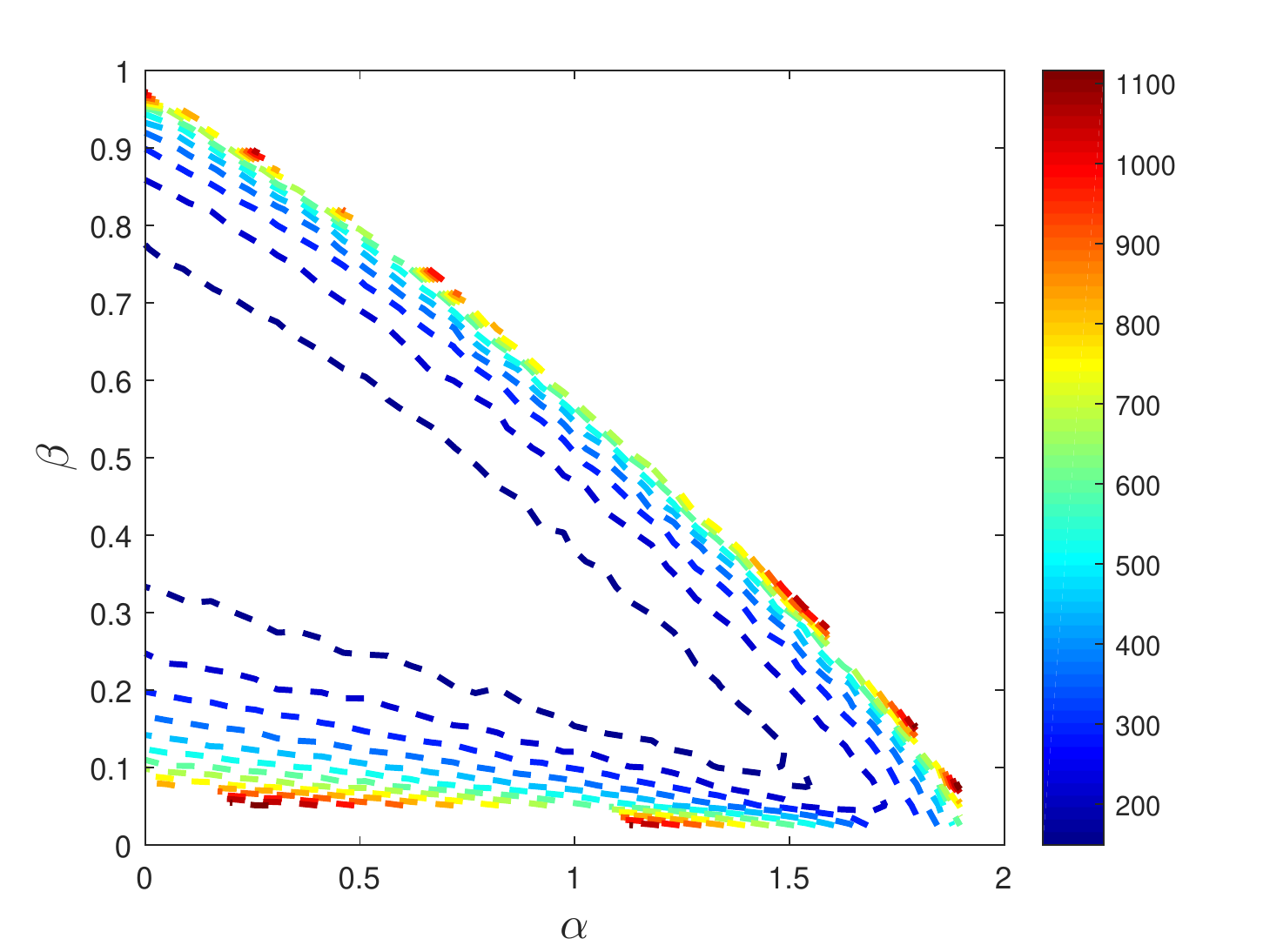}}
	\subfigure[mFDBK: $n=50$ and $m=50n$]{
	    \includegraphics[width=0.48\textwidth]{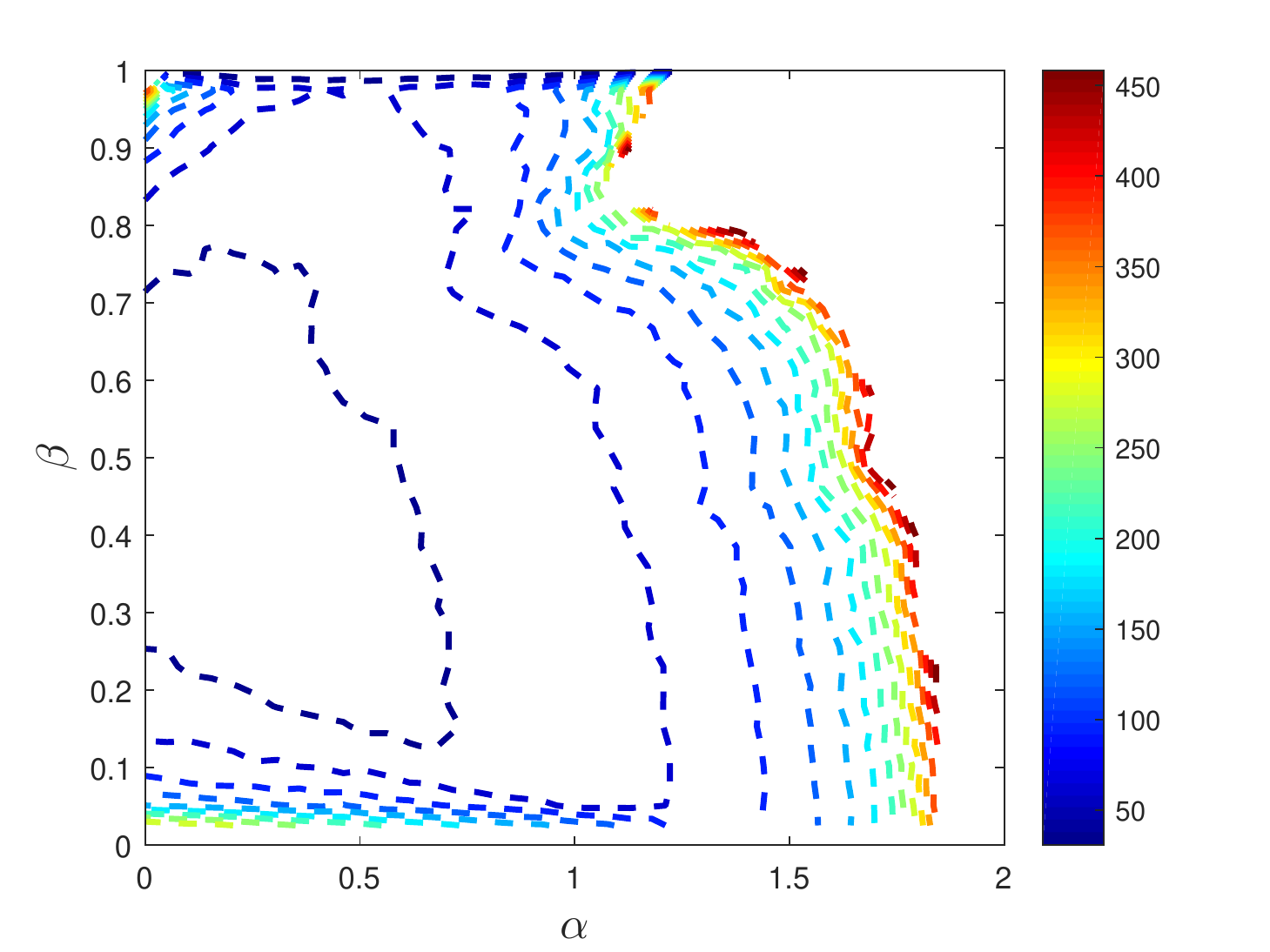}}
    \subfigure[mMWRK: $n=50$ and $m=100n$]{
		\includegraphics[width=0.48\textwidth]{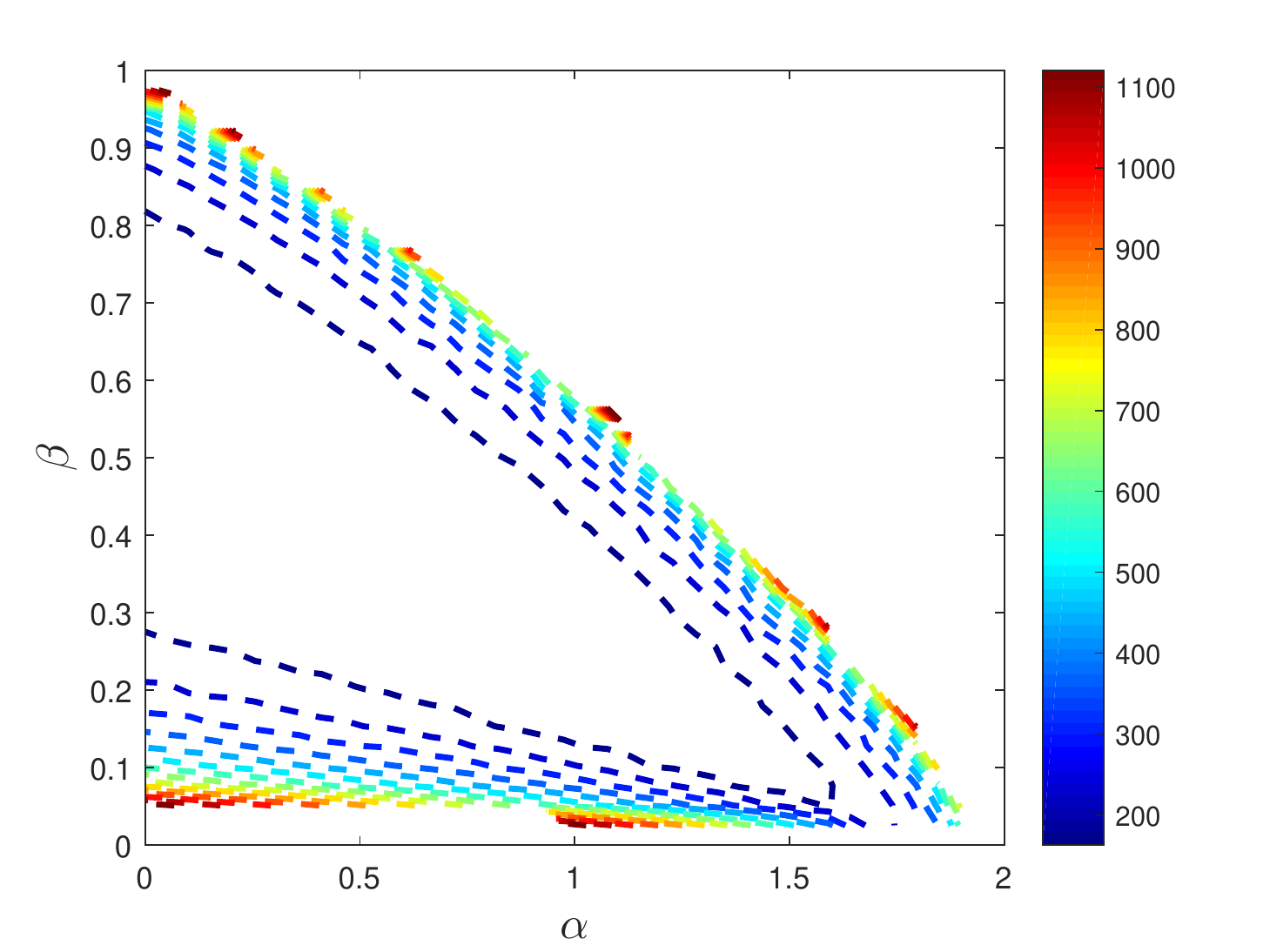}}
	\subfigure[mFDBK: $n=50$ and $m=100n$]{
	    \includegraphics[width=0.48\textwidth]{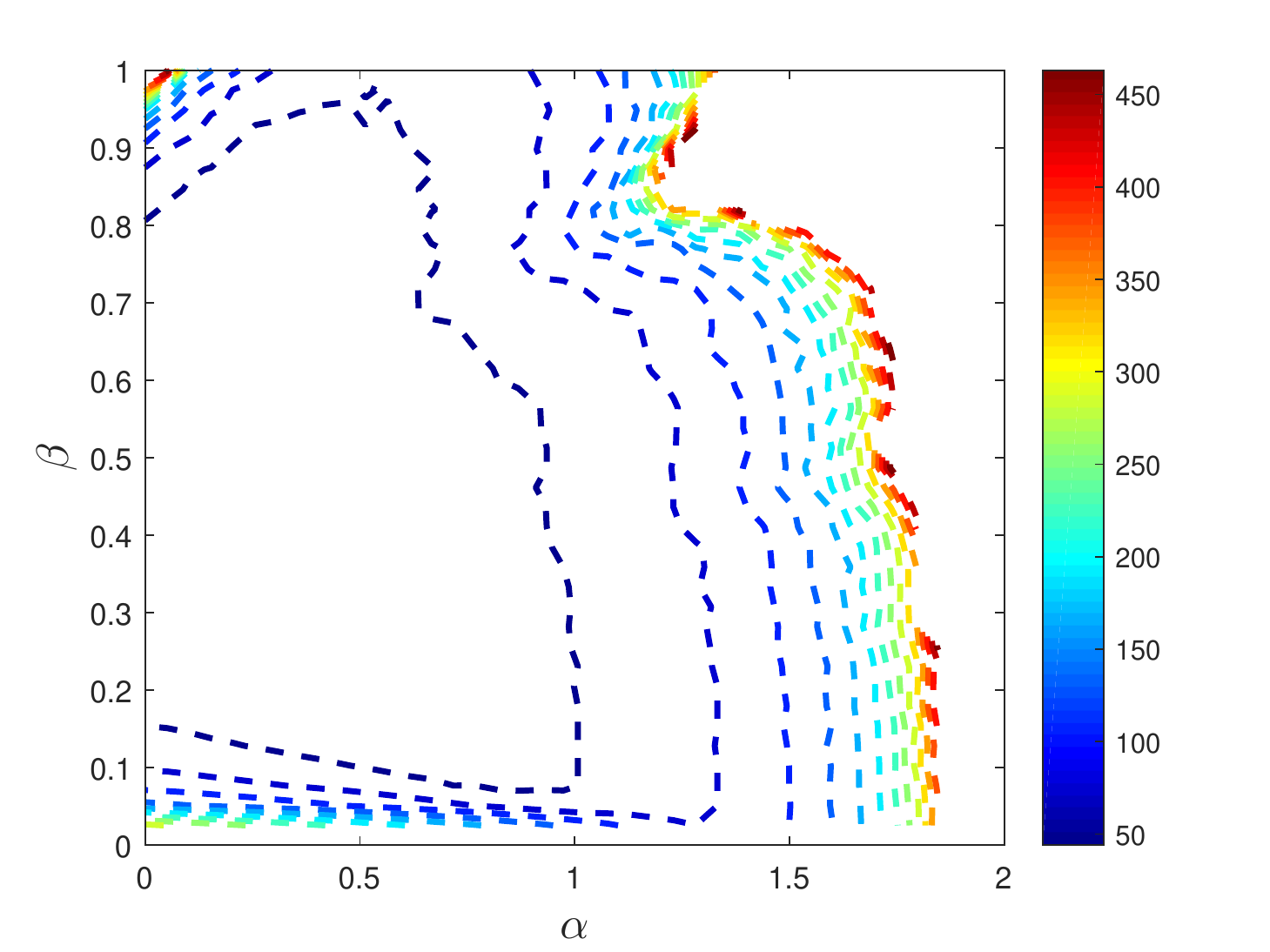}}
\caption{The number of iteration steps of mMWRK (left) and mFDBK (right) with different $(\alpha,\beta)$ for solving a consistent linear systems, where the test matrix is generated by {\tt randn}$(m,n)$.}
\label{ContourFig_Alpha_Beta_VS_IT_randn_m_Larger}
\end{figure}

\begin{figure}[!htb]
\centering
    \subfigure[mMWRK: $m=50$ and $n=10m$]{
		\includegraphics[width=0.48\textwidth]{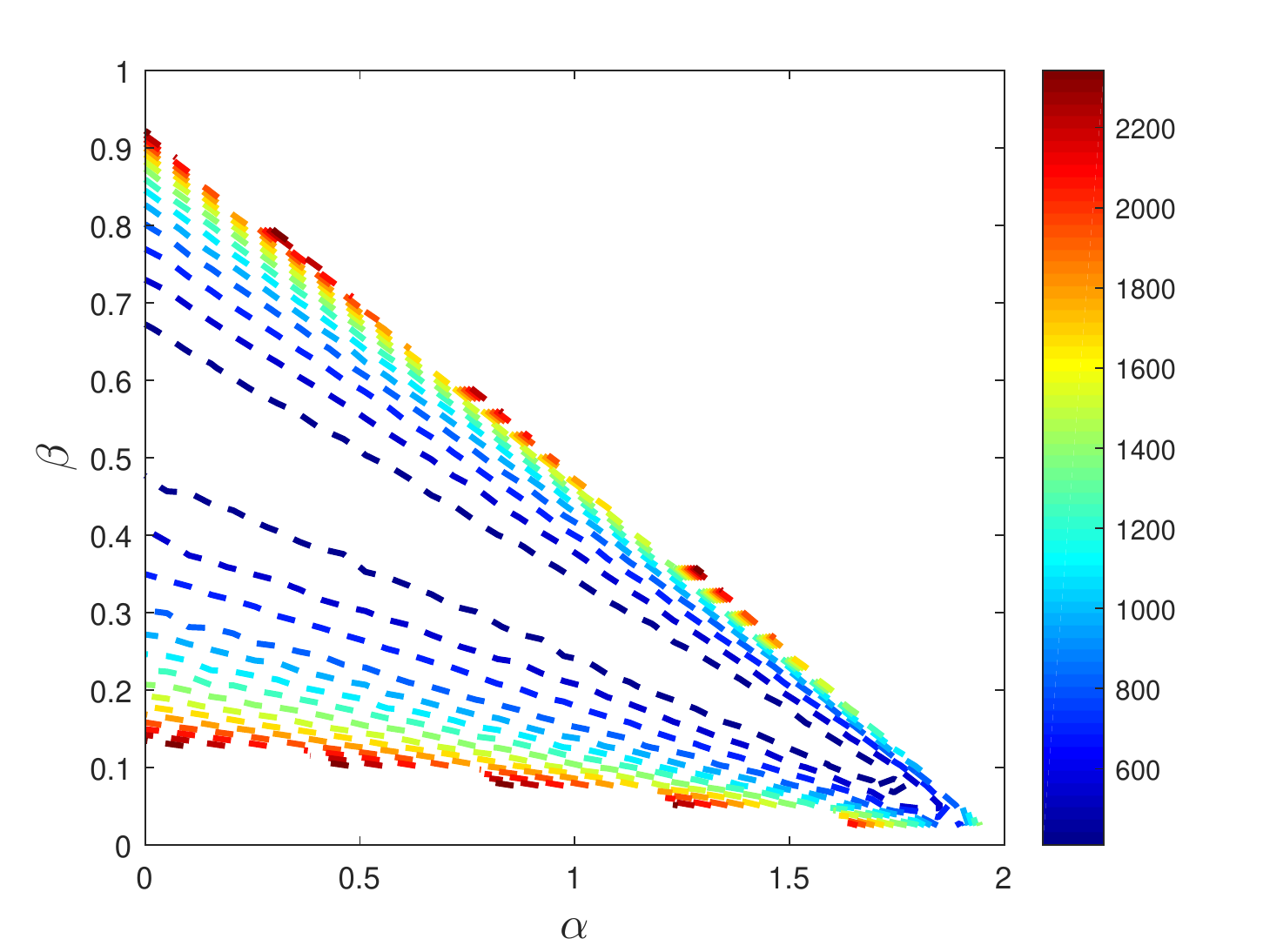}}
	\subfigure[mFDBK: $m=50$ and $n=10m$]{
	    \includegraphics[width=0.48\textwidth]{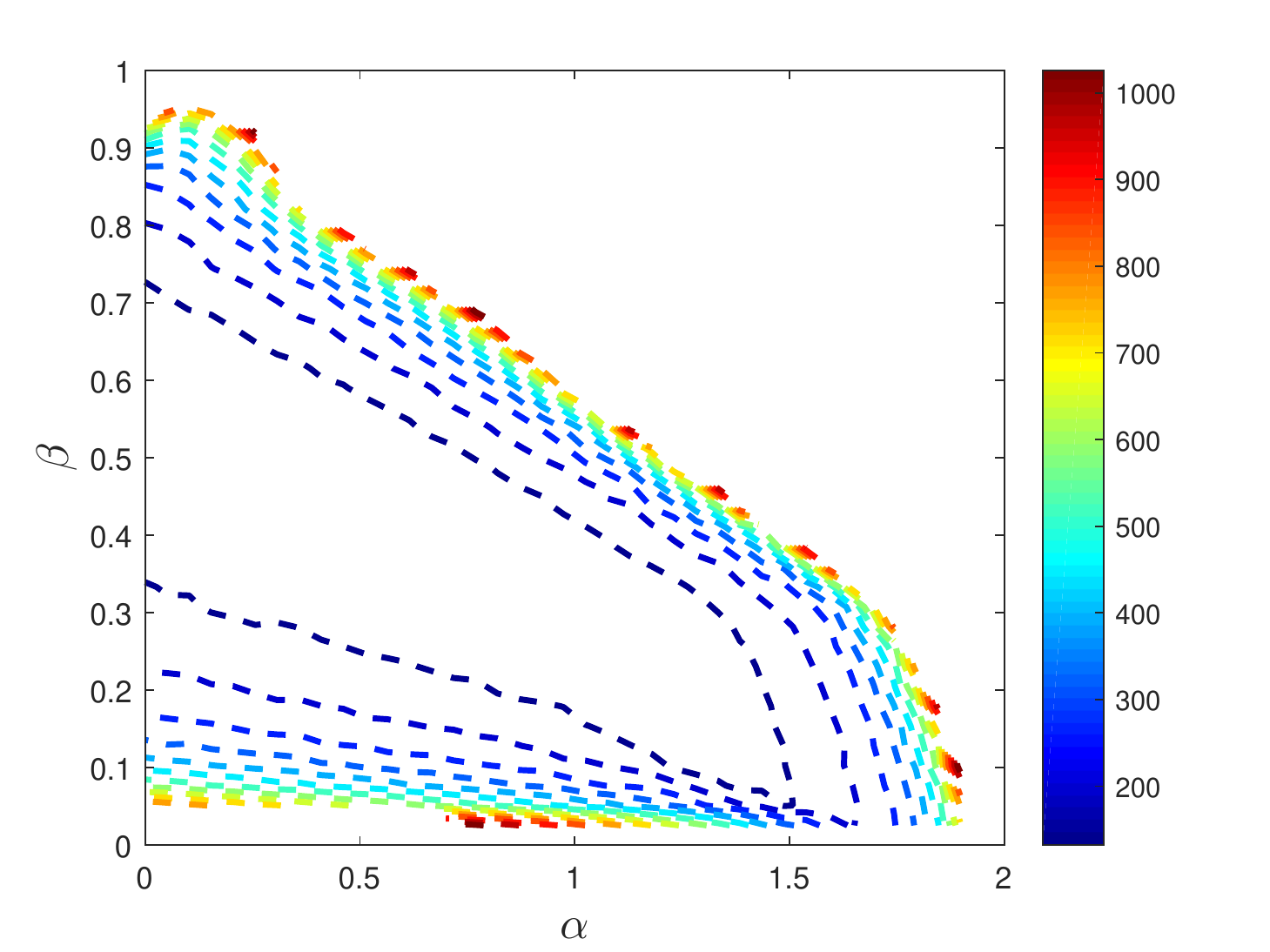}}
    \subfigure[mMWRK: $m=50$ and $n=50m$]{
		\includegraphics[width=0.48\textwidth]{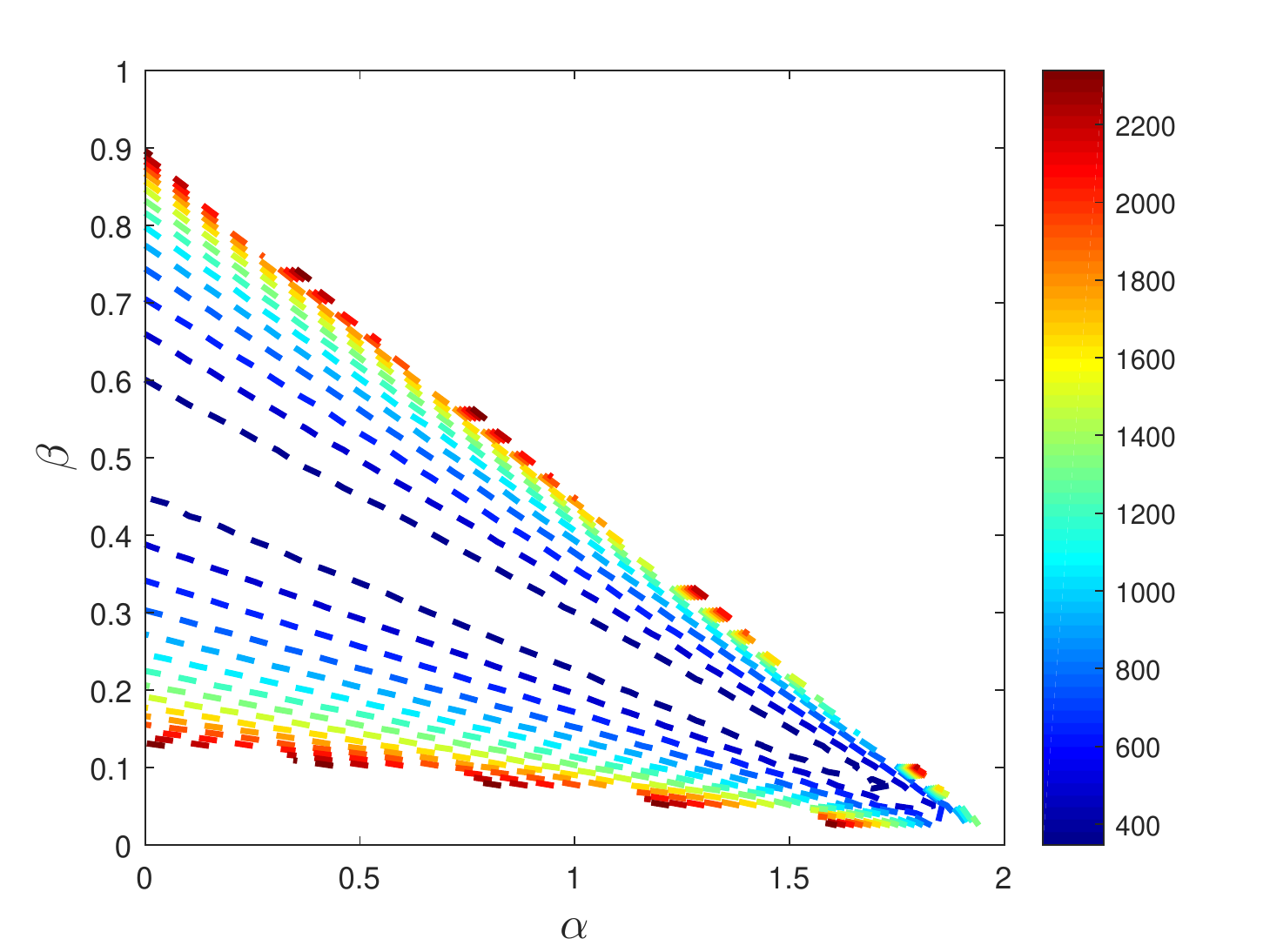}}
	\subfigure[mFDBK: $m=50$ and $n=50m$]{
	    \includegraphics[width=0.48\textwidth]{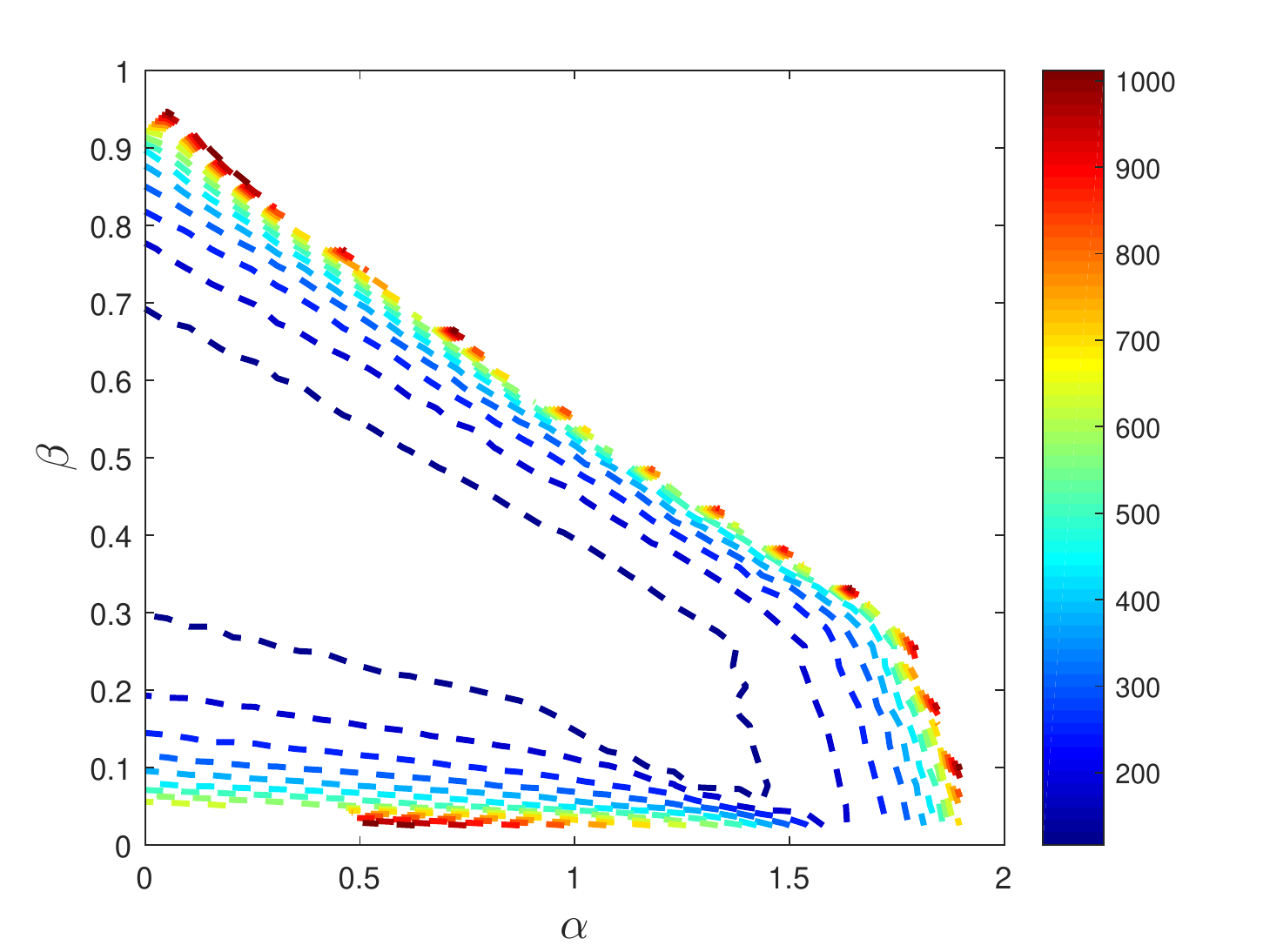}}
    \subfigure[mMWRK: $m=50$ and $n=100m$]{
		\includegraphics[width=0.48\textwidth]{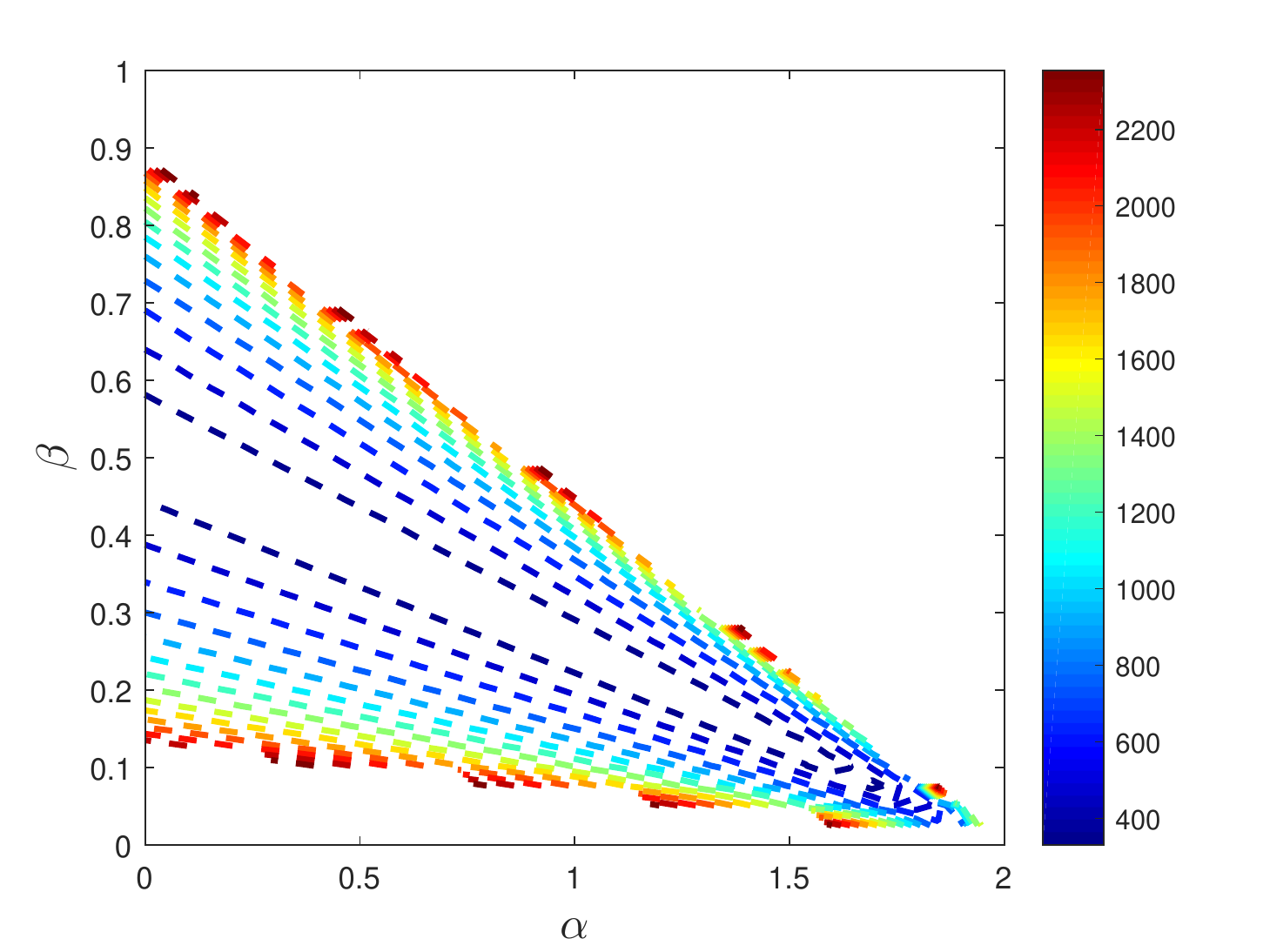}}
	\subfigure[mFDBK: $m=50$ and $n=100m$]{
	    \includegraphics[width=0.48\textwidth]{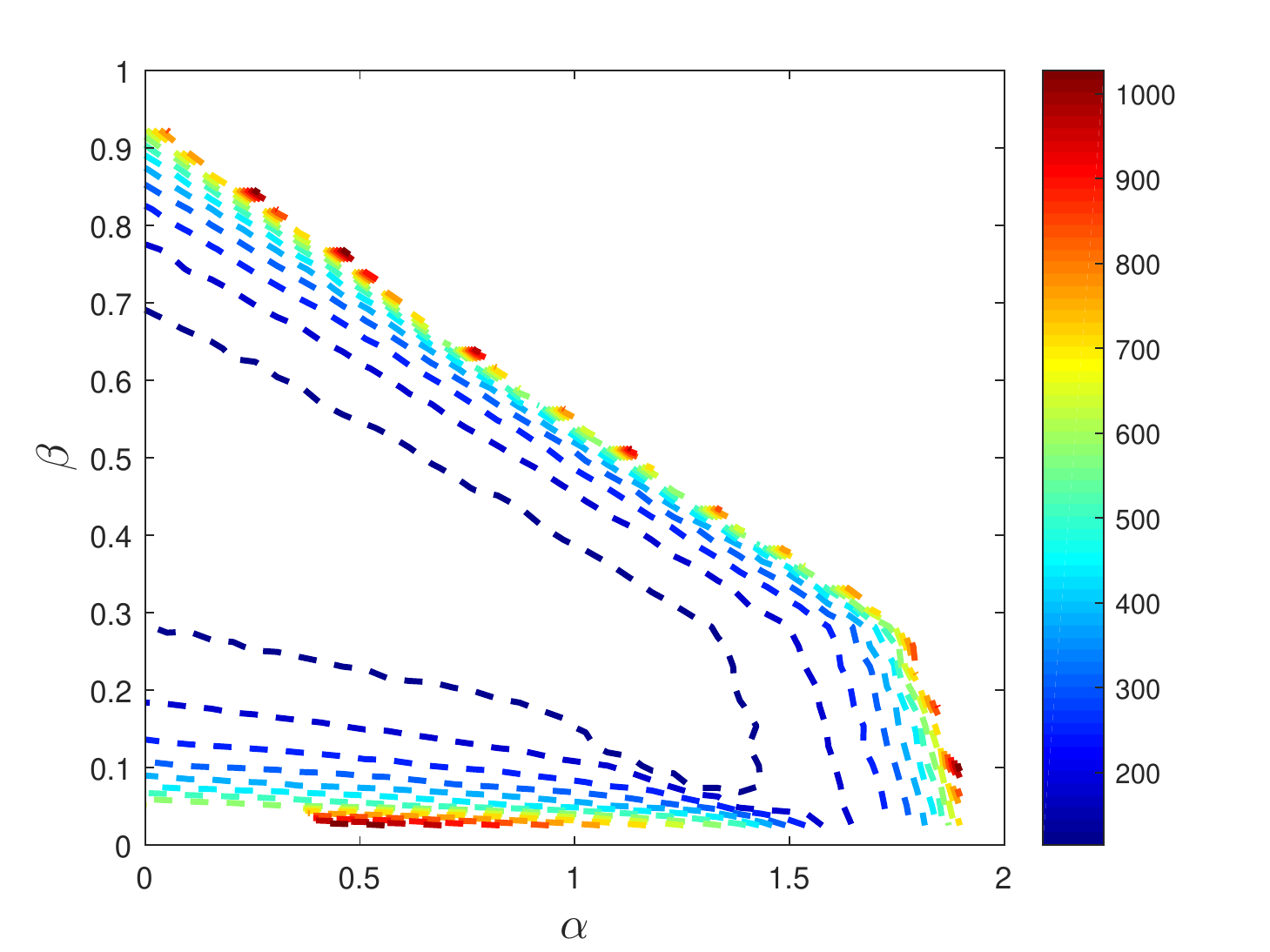}}
\caption{The number of iteration steps of mMWRK (left) and mFDBK (right) with different $(\alpha,\beta)$ for solving a consistent linear systems, where the test matrix is generated by {\tt randn}$(m,n)$.}
\label{ContourFig_Alpha_Beta_VS_IT_randn_n_Larger}
\end{figure}

To begin with, we depict the performances, given by mMWRK (left) and mFDBK (right) solving the over-determine linear systems, in Figure \ref{ContourFig_Alpha_Beta_VS_IT_randn_m_Larger} with $n=50$ when $m=10n$, $50n$, and $100n$. In this figure, the number of iteration steps is represented by the contour lines. Furthermore, we test our algorithm on solving under-determined linear systems. To do so, we set $m=50$ with $n=10m$, $50m$, and $100m$. The results are shown in Figure \ref{ContourFig_Alpha_Beta_VS_IT_randn_n_Larger}. In these two figures, we can observe the following phenomena. (I) The momentum method needs less number of iteration steps than the corresponding original no-momentum variant by setting  $(\alpha,\beta)=(1,0)$. It implies that the momentum technique can further improve the convergence behavior of the greedy deterministic single and multiple row-action methods. (II) The results show that the parameter pairs $(\alpha,\beta)=(0.75,0.5)$ and $(0.5,0.5)$ are the good choices for mMWRK and mFDBK, respectively, and result in a satisfactory convergence. Not especially specified, we will adopt this parameter
 selection approach for mMWRK and mFDBK in the following numerical test.

\subsection{Synthetic data} The following coefficient matrix is generated from synthetic data, which is dense and can yield various specific instantiations about the linear systems  (full rank or rank-deficient, over- or under-determined) by varying the input parameters.

\begin{example}\label{ERMR:example1}
  As in Du et al. \cite{20DSS}, for given $m$, $n$, $r$, and $\kappa>1$, we construct a dense matrix $A$ by $A=UDV^T$, where $U\in \mathbb{R}^{m\times r}$, $D\in \mathbb{R}^{r\times r}$, and $V\in \mathbb{R}^{n\times r}$. Using MATLAB colon notation, these matrices are generated by
  $[U,\sim]={\rm qr}({\rm randn}(m,r),0)$,
  $[V,\sim]= {\rm qr}({\rm randn}(n,r),0)$, and
  $D={\rm diag}(1 + (\kappa - 1).* {\rm rand}(r,1))$.
\end{example}

In the following, we consider two types of {\it rank-deficient} cases by setting (I) $m>n$, $r=n/2$, and $\kappa=n/10$; (II) $m<n$, $r=m/2$, and $\kappa=m/10$. For solving the linear system \eqref{eq:Ax=b} concerning this class of coefficient matrices,
we list the number of iteration steps in Tables \ref{tab:Ex_randn_m_larger_n}-\ref{tab:Ex_randn_n_larger_m}. Since the coefficient matrix given by Example \ref{ERMR:example1} has randomness, we repeat $20$ runs of the MWRK, mMWRK, FDBK, and mFDBK methods, and present the median of the results. The two tables reveal that the mMWRK (resp., mFDBK) method has a faster convergence rate than the MWRK (resp., FDBK) method. The maximum of the speed-ups is $1.63$ (resp., $2.37$) and the minimum is $1.50$ (resp., $1.29$).

\begin{table}[!htb]
\caption{The numerical results obtained by MWRK, mMWRK, FDBK, and mFDBK  to solve the linear system \eqref{eq:Ax=b} in Example \ref{ERMR:example1} for $m > n$, $r = n/10$, and $\kappa = n/10$.}
\centering
\begin{tabular}{ccccc}
\cline{1-5}
$m\times n$&
$10000 \times 350$&$10000 \times 400$&
$10000 \times 450$&$10000 \times 500$\\
\cline{1-5}
MWRK  &4146.8&6240.2&9496.6&11565.0\\
mMWRK &2771.2&4103.0&6235.8&7577.4\\
SU$_1$&{\bf 1.50}&{\bf 1.52}&{\bf 1.52}&{\bf 1.53}\\
\cline{1-5}
FDBK  &2029.0&2618.8&3848.0&4206.8\\
mFDBK &1068.2&1430.8&1850.0&2115.2\\
SU$_2$&{\bf 1.90}&{\bf 1.83}&{\bf 2.08}&{\bf 1.99}\\
\cline{1-5}
$m\times n$&
$15000 \times 350$&$15000 \times 400$&
$15000 \times 450$&$15000 \times 500$\\
\cline{1-5}
MWRK  &4045.1&4938.8&7226.6&9466.8\\
mMWRK &2647.5&3235.8&4754.4&6239.8\\
SU$_1$&{\bf 1.53}&{\bf 1.53}&{\bf 1.52}&{\bf 1.52}\\
\cline{1-5}
FDBK  &1735.2&2119.6&3055.6&3384.0\\
mFDBK &732.0&1026.4&1589.2&1744.6\\
SU$_2$&{\bf 2.37}&{\bf 2.07}&{\bf 1.92}&{\bf 1.94}\\
\cline{1-5}
$m\times n$&
$20000 \times 350$&$20000 \times 400$&
$20000 \times 450$&$20000 \times 500$\\
\cline{1-5}
MWRK   &3191.4&4766.4&6371.4&8846.8\\
mMWRK  &2076.6&3127.6&4213.4&5785.2\\
SU$_1$ &{\bf 1.54}&{\bf 1.52}&{\bf 1.51}&{\bf 1.53}\\
\cline{1-5}
FDBK   &1655.6&2231.4&2658.6&3340.6\\
mFDBK  &809.8&1028.8&1330.8&1535.6\\
SU$_2$ &{\bf 2.04}&{\bf 2.17}&{\bf 2.00}&{\bf 2.18}\\
\cline{1-5}
\end{tabular}
\label{tab:Ex_randn_m_larger_n}
\end{table}

\begin{table}[!htb]
\caption{The numerical results obtained by MWRK, mMWRK, FDBK, and mFDBK to solve the linear system \eqref{eq:Ax=b} in Example \ref{ERMR:example1} for $n>m$, $r = m/10$, and $\kappa = m/10$.}
\centering
\begin{tabular}{ccccc}
\cline{1-5}
$m\times n$&
$350 \times 10000$&$400 \times 10000$&
$450 \times 10000$&$500 \times 10000$\\
\cline{1-5}
MWRK   &10853.8&14463.1&19904.0&30162.4\\
mMWRK  &6775.4&8990.0&12364.1&18468.1\\
SU$_1$ &{\bf 1.60}&{\bf 1.61}&{\bf 1.61}&{\bf 1.63}\\
\cline{1-5}
FDBK   &3258.2&4451.5&5706.3&7978.9\\
mFDBK  &2523.0&3368.5&4207.5&6041.1\\
SU$_2$ &{\bf 1.29}&{\bf 1.32}&{\bf 1.36}&{\bf 1.32}\\
\cline{1-5}
$m\times n$&
$350 \times 15000$&$400 \times 15000$&
$450 \times 15000$&$500 \times 15000$\\
\cline{1-5}
MWRK   &10341.9&9874.6&19081.0&29671.4\\
mMWRK  &6475.8&6104.8&11809.4&18367.1\\
SU$_1$ &{\bf 1.60}&{\bf 1.62}&{\bf 1.62}&{\bf 1.62}\\
\cline{1-5}
FDBK   &2491.1&3087.1&5255.2&7721.5\\
mFDBK  &1518.7&2346.6&4000.6&5849.5\\
SU$_2$ &{\bf 1.64}&{\bf 1.32}&{\bf 1.31}&{\bf 1.32}\\
\cline{1-5}
$m\times n$&
$350 \times 20000$&$400 \times 20000$&
$450 \times 20000$&$500 \times 20000$\\
\cline{1-5}
MWRK   &4599.9&12399.9&12935.0&13773.6\\
mMWRK  &2880.3&7686.1&8054.4&8430.4\\
SU$_1$ &{\bf 1.60}&{\bf 1.61}&{\bf 1.61}&{\bf 1.63}\\
\cline{1-5}
FDBK   &1500.6&3834.4&3793.1&3526.4\\
mFDBK  &1151.6&2845.6&2855.5&2672.9\\
SU$_2$ &{\bf 1.30}&{\bf 1.35}&{\bf 1.33}&{\bf 1.32}\\
\cline{1-5}
\end{tabular}
\label{tab:Ex_randn_n_larger_m}
\end{table}

\begin{figure}[!htb]
\centering
    \subfigure[mMWRK]{
		\includegraphics[width=0.48\textwidth]{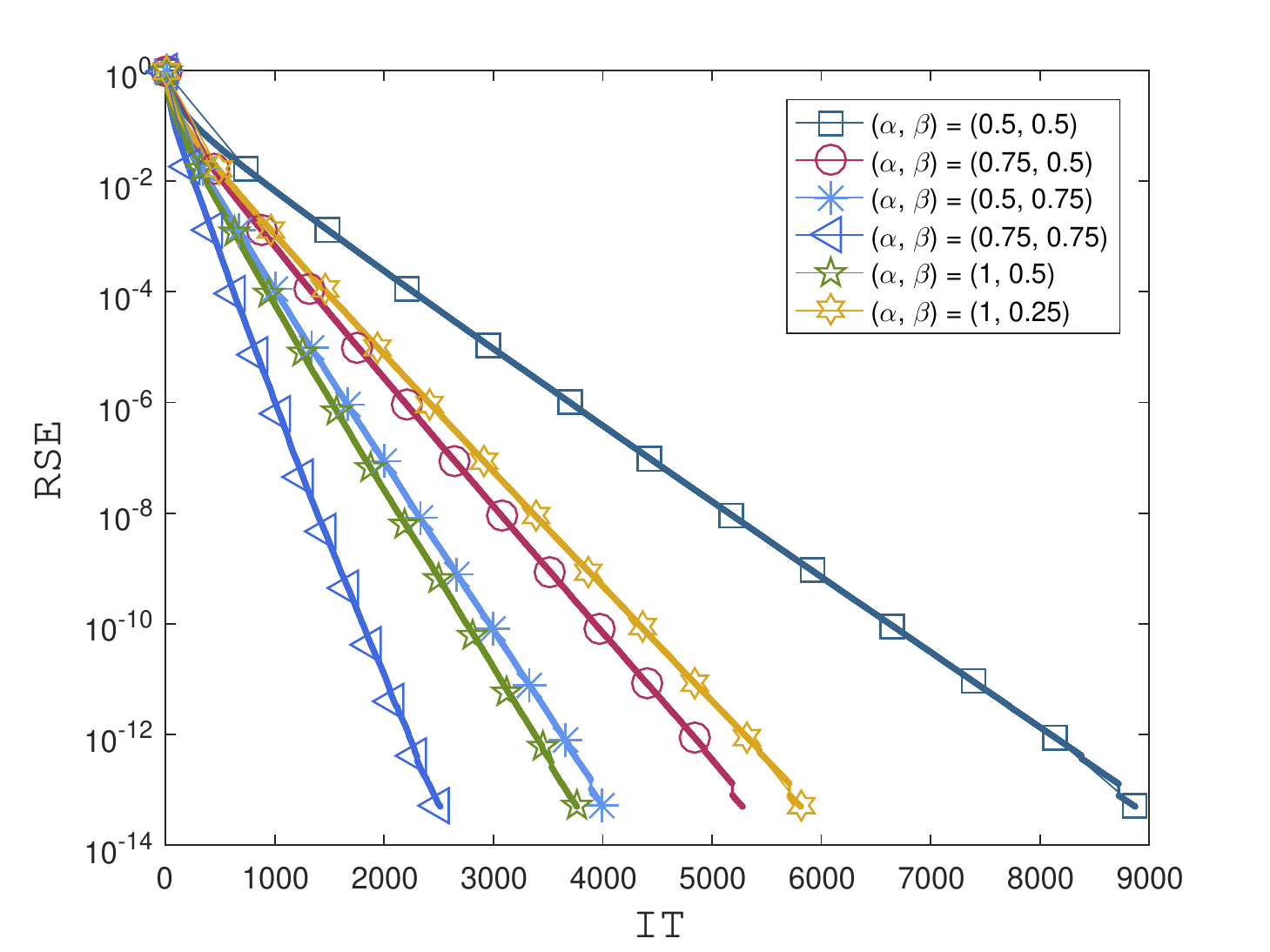}}
	\subfigure[mFDBK]{
	    \includegraphics[width=0.48\textwidth]{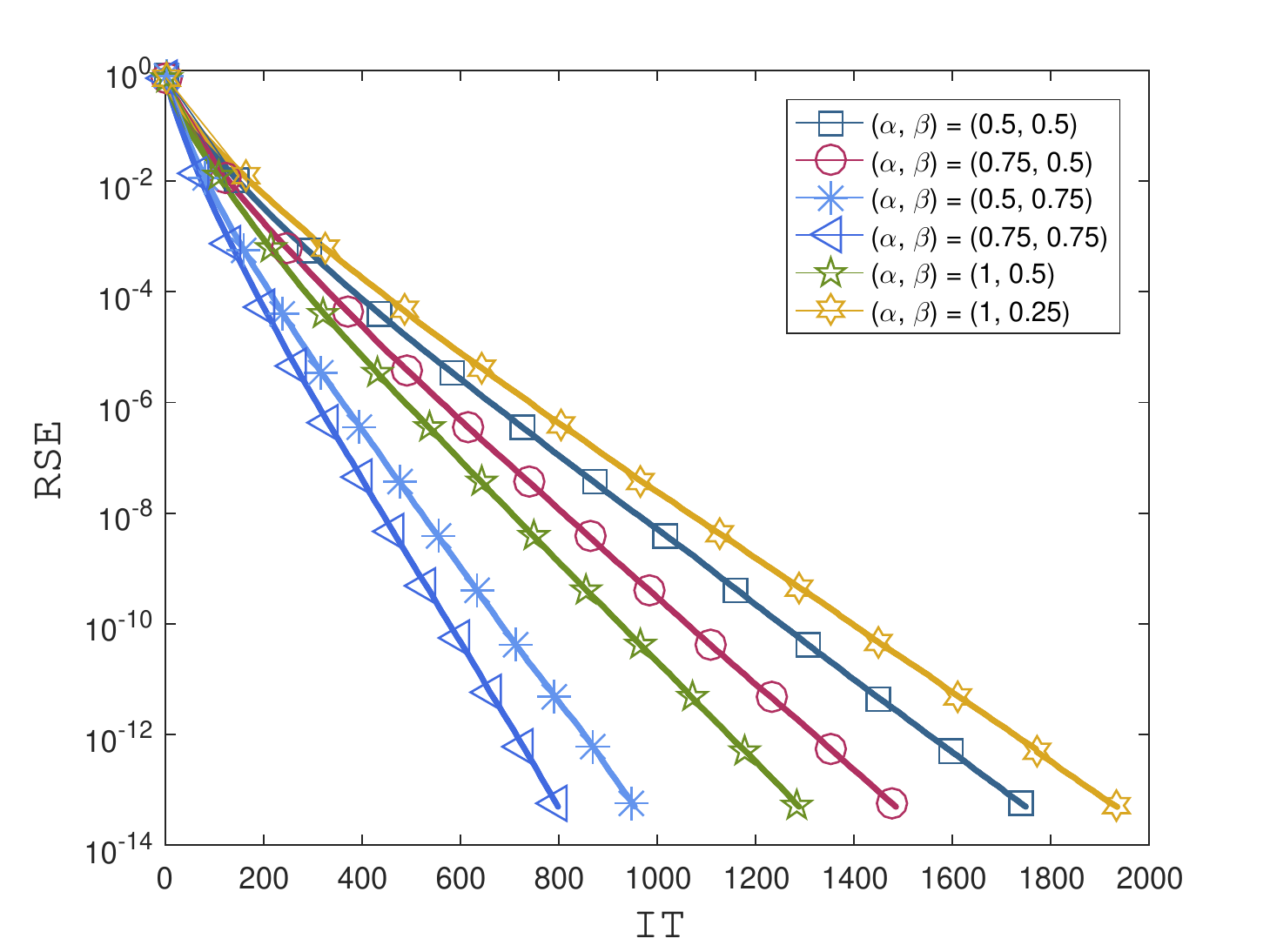}}
\caption{RSE versus IT obtained by  mMWRK (a) and mFDBK (b) for Example \ref{ERMR:example1} when $m=15000$, $n=350$, $r = n/10$, and $\kappa = n/10$.}
\label{mMWRK_m_n350_variousAlpha+Beta}
\end{figure}

We emphasize that the experimentally iterative parameter pair $(\alpha,\beta)$, given by Section \ref{sec:choice_alpha+beta}, are used and this selection is not optimal. We may choose them to be any positive constant bounded by Remark \ref{Remark:gamma1+gamma2<1}. Selecting an appropriate pair of iteration parameters may allow the momentum method to converge more quickly. For the sake of illustration, we utilize the mMWRK and mFDBK methods with various $(\alpha,\beta)$ to solve an over-determined  linear system \eqref{eq:Ax=b} in Example \ref{ERMR:example1} and depict their convergence behaviors of RSE versus IT in Figure \ref{mMWRK_m_n350_variousAlpha+Beta}. In this example, we assign values to the input parameters as $m=15000$, $n=350$, $r = n/10$, and $\kappa = n/10$, and to $(\alpha,\beta)$ as $(0.5,0.5)$, $(0.75,0.5)$, $(0.5,0.75)$, $(0.75,0.75)$, $(1,0.25)$, and $(1,0.5)$. It can be seen that both mMWRK and mFDBK successfully compute an approximate solution for all cases. The fastest convergence rate occurs if one takes the parameter pair $(\alpha,\beta) = (0.75, 0.75)$. In this case, the iteration counts of mMWRK and mFDBK are $1225.8$ and $389.6$, respectively, and the corresponding speed-ups are SU$_1=3.30$ and SU$_2=4.45$, which are appreciably larger than those in Table \ref{tab:Ex_randn_m_larger_n}.

\begin{figure}[!htb]
\centering
    \subfigure[mMWRK]{
		\includegraphics[width=0.48\textwidth]{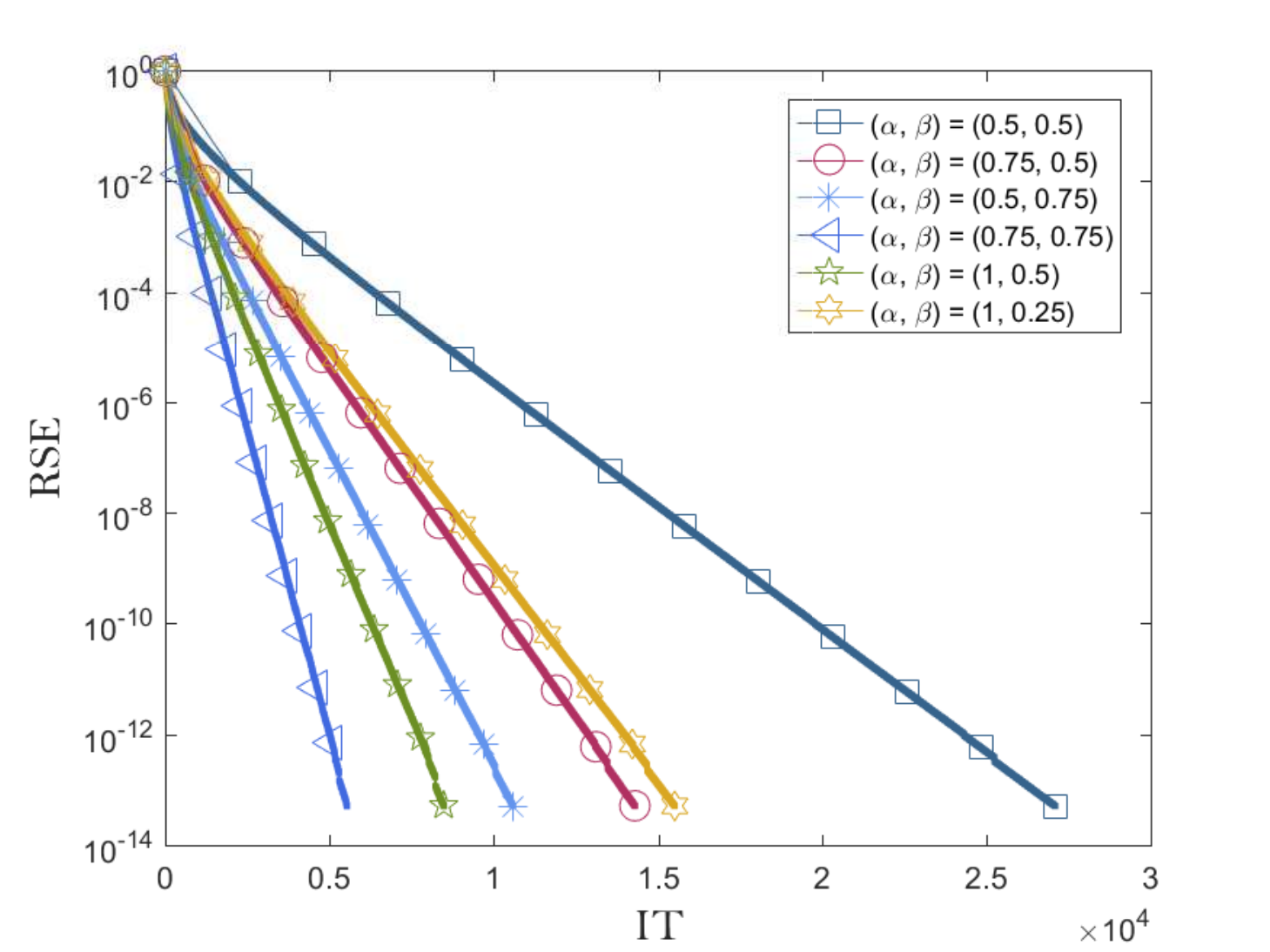}}
	\subfigure[mFDBK]{
	    \includegraphics[width=0.48\textwidth]{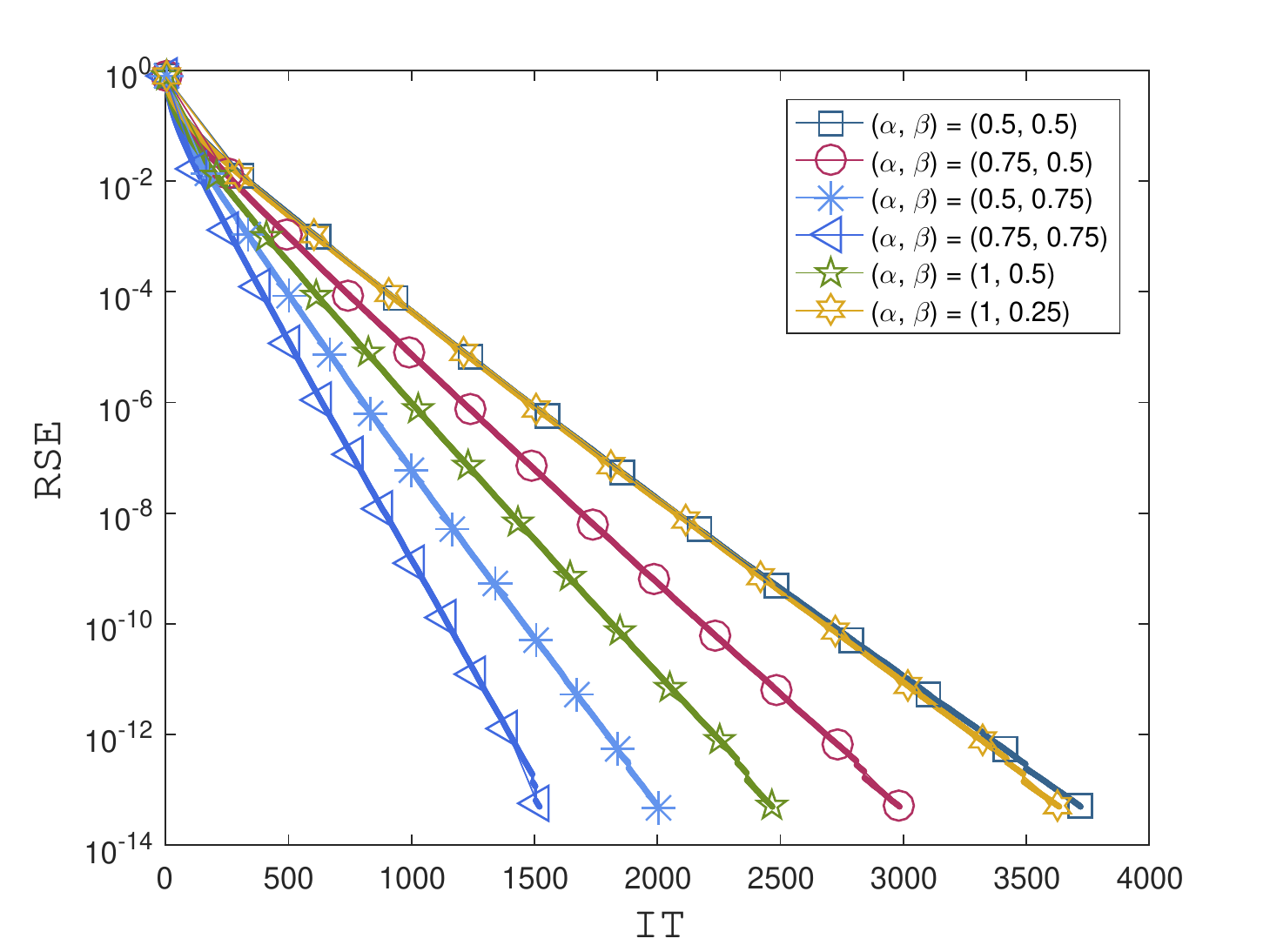}}
\caption{RSE versus IT obtained by  mMWRK (a) and mFDBK (b) for Example \ref{ERMR:example1} when $m=350$, $n=15000$, $r = m/10$, and $\kappa = m/10$.}
\label{mMWRK_m350_n_variousAlpha+Beta}
\end{figure}

By setting the input parameter as $m=350$, $n=15000$, $r = m/10$, and $\kappa = m/10$ for Example \ref{ERMR:example1}, we obtain an under-determined linear system. We choose the same parameter pairs as above for mMWRK and mFDBK, and plot the corresponding convergence behavior of RSE versus IT in Figure \ref{mMWRK_m350_n_variousAlpha+Beta}. As the figure depicts, the numerical phenomena are similar to those described above.  That is,  mMWRK and mFDBK arrive at the fastest convergence rate with $(\alpha,\beta) = (0.75, 0.75)$. For this case, the iteration counts of mMWRK and mFDBK are considerably smaller than those in Table \ref{tab:Ex_randn_n_larger_m}, with the speed-ups being $3.84$ (IT $=2692.4$) and $3.18$ (IT $=783.2$), respectively.

\subsection{Real-world application: curve fitting} The subsequent numerical experiments consider the linear systems in computer-aided geometric design (CAGD), such as curve fitting.

Let us fit the ordered point set $ \{q_i \in \mathbb{R}^3: i\in[m] \}$.  Assume that $\left\{\mu_j(x):j\in [n]\right\}$
is a basis sequence and $\{p_j^{(k)} \in \mathbb{R}^3: j \in[n] \}$ is the control point sequence at $k$th iterate. The $k$th fitting curve, generated by the geometric iterative method (GIM) in \cite{18LMD},
\begin{align}\label{eq:curve_C(x)}
  \mathrm{C}^{(k)}(x) =\sum_{j=1}^n \mu_j(x) p_j^{(k)},~x \in [x_1,~x_m],
\end{align}
progressively approximates a target curve by updating the control points according to
\begin{align*}
 p_j^{(k+1)} = p_j^{(k)} + \delta_j^{(k)},
\end{align*}
where $\delta_j^{(k)} \in \mathbb{R}^3$ is called the adjust vector and computed by $q_i$ and $p_j^{(k)}$.

Let the $x$-, $y$-, and $z$-coordinates of control point $p_j^{(k)}$ (resp., adjust vector $\delta_j^{(k)}$) be respectively stored in the vectors $p_x^{(k)}$, $p_y^{(k)}$, and $p_z^{(k)}$  (resp., $\delta_x^{(k)}$, $\delta_y^{(k)}$, and $\delta_z^{(k)}$). From algebraic aspects, the GIM iterative processes,
\begin{align*}\renewcommand\arraystretch{1.5}
   p_x^{(k+1)} = p_x^{(k)} + \delta_x^{(k)},~
   p_y^{(k+1)} = p_y^{(k)} + \delta_y^{(k)},~{\rm and}~~
   p_z^{(k+1)} = p_z^{(k)} + \delta_z^{(k)},
\end{align*}
are equal to iteratively solving three linear systems. Therefore, the mMWRK and mFDBK method are suitable for addressing this issue.

The implementation detail of mMWRK and mFDBK curve fittings are presented as follows. Let the data points be arranged into
$q = \left[q_{1}~q_{2}~\cdots ~q_{m}\right]^T =\left[ q_x ~ q_y ~ q_z\right] \in \mathbb{R}^{m\times 3}$. We input
the collocation matrix $A$,
two initial vectors $p_x^{(1)}=p_x^{(0)}\in \mathbb{R}^{n}$ (resp., $p_y^{(1)}=p_y^{(0)}$, $p_z^{(1)}=p_z^{(0)}$),
the right-hand side $q_x\in \mathbb{R}^{m}$ (resp., $q_y$, $q_z$),
and compute the next vector $p_x^{(k)}$ (resp., $p_y^{(k)}$, $p_z^{(k)}$) using the mMWRK and mFDBK update rules. Then, the approximate curve is formulated according to formula \eqref{eq:curve_C(x)}.

\begin{example}\label{mMWRK and mFDBK:example2}
We fit the data points $ \{q_i: i\in[m] \}$ sampled from the following curves,
\begin{equation*}
\begin{array}{ll}
{\rm Curve~1:}
&x = 30\cos( t\pi/3), \vspace{1ex}\\
&y = 30\sin(t\pi/3),  \vspace{1ex}\\
&z = 3 t\pi~~(0 \leq t \leq 10\pi); \vspace{1ex}\\
{\rm Curve~2:}
&x = -22 \cos(t) - 128 \sin(t) - 44 \cos(3t) - 78 \sin(3t), \vspace{1ex}\\
&y = -10 \cos(2t) - 27 \sin(2t) + 38 \cos(4t) + 46 \sin(4t),  \vspace{1ex}\\
&z =  70 \cos(3t) - 40 \sin(3t) ~~ (0 \leq t \leq 2 \pi),
\end{array}
\end{equation*}
which come from reference \cite{14DL} and the collection of various topics in geometry (available from \verb"http://paulbourke.net/geometry/"). As other researchers do, we first assign a parameter sequence $\nu$ and a knot vector $\mu$ of cubic B-spline basis, which is simple and has a wide range of applications in CAGD; see, e.g., \cite{18LMD}, and then obtain the collocation matrix by using the MATLAB built-in function, e.g., $A={\rm spcol}(\mu, 4, \nu)$.
\end{example}

\begin{figure}[!htb]
\centering
    \subfigure[Curve 1]{
		\includegraphics[width=0.48\textwidth]{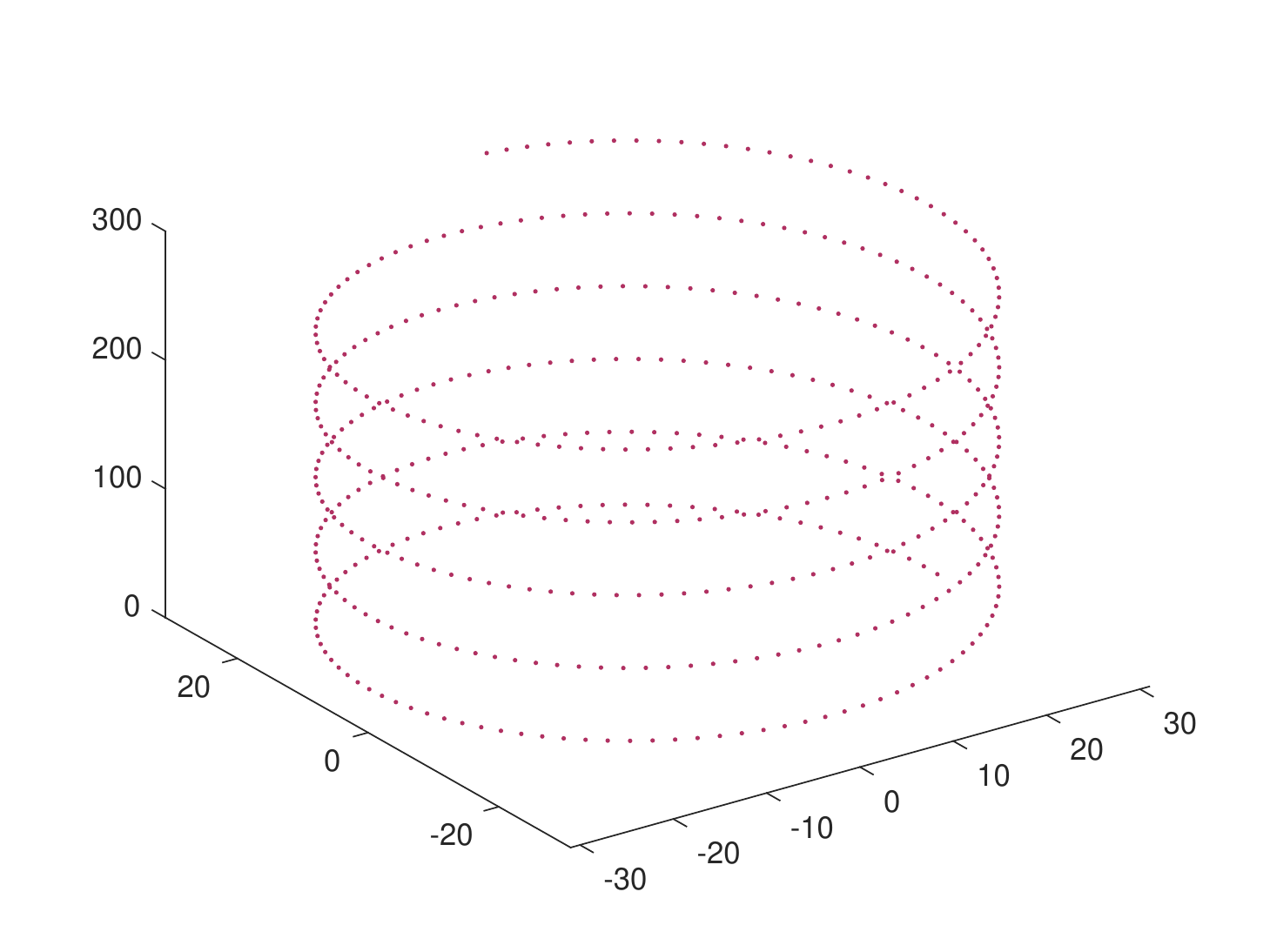}}
	\subfigure[Curve 2]{
	    \includegraphics[width=0.48\textwidth]{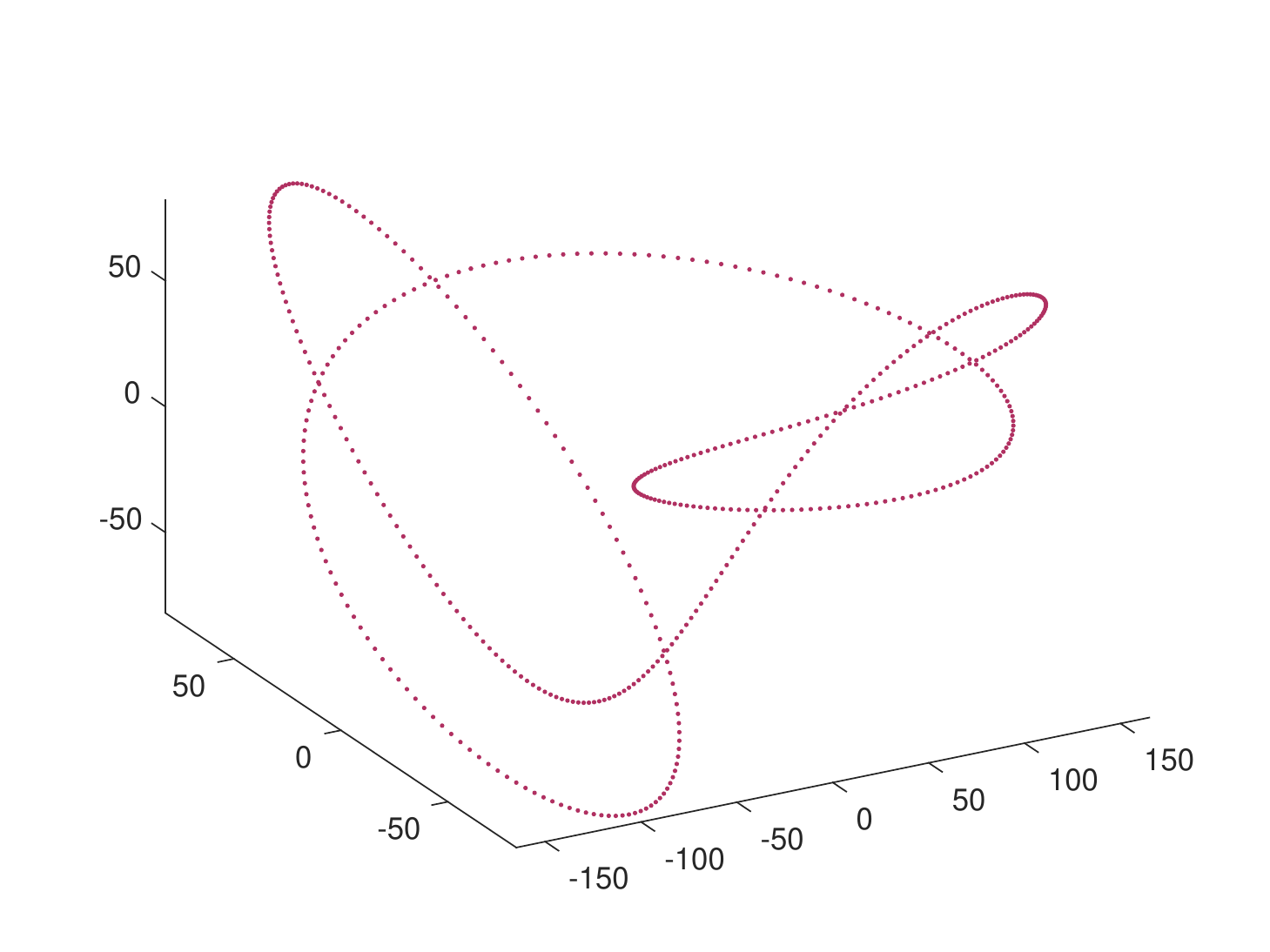}}
\caption{$500$ initial data points sampled from Curve 1 (a) and 2 (b) in Example \ref{mMWRK and mFDBK:example2}.}
\label{fig:granny knot_curve}
\end{figure}

The mMWRK and mFDBK methods can be started with arbitrary initial control points in the column space of the collocation matrix, and a suitable and efficient choice is to set the initial control points to be zeros.  At $k$th iterate, the relative solution error is defined by RSE = $E_k/E_0$ for $k=0,1,2,\cdots$, where $E_k = \bF{ p^{(k)} -  p_{\ast}}$ and $p_{\ast} = A^{\dag} q$ is the least-square solution. The computation  is terminated once RSE is less than $10^{-12}$.

The initial data points sampled from Curves 1-2 in Example \ref{mMWRK and mFDBK:example2} for $m = 500$ are shown in Figure \ref{fig:granny knot_curve} as a concrete example. In the following, $n$ control points are employed and the capabilities of mMWRK and mFDBK to fit  $m$ three-dimensional data points are discussed.

\begin{figure}[!htb]
\centering
	\subfigure[The mMWRK curve]{
	    \includegraphics[width=0.48\textwidth]{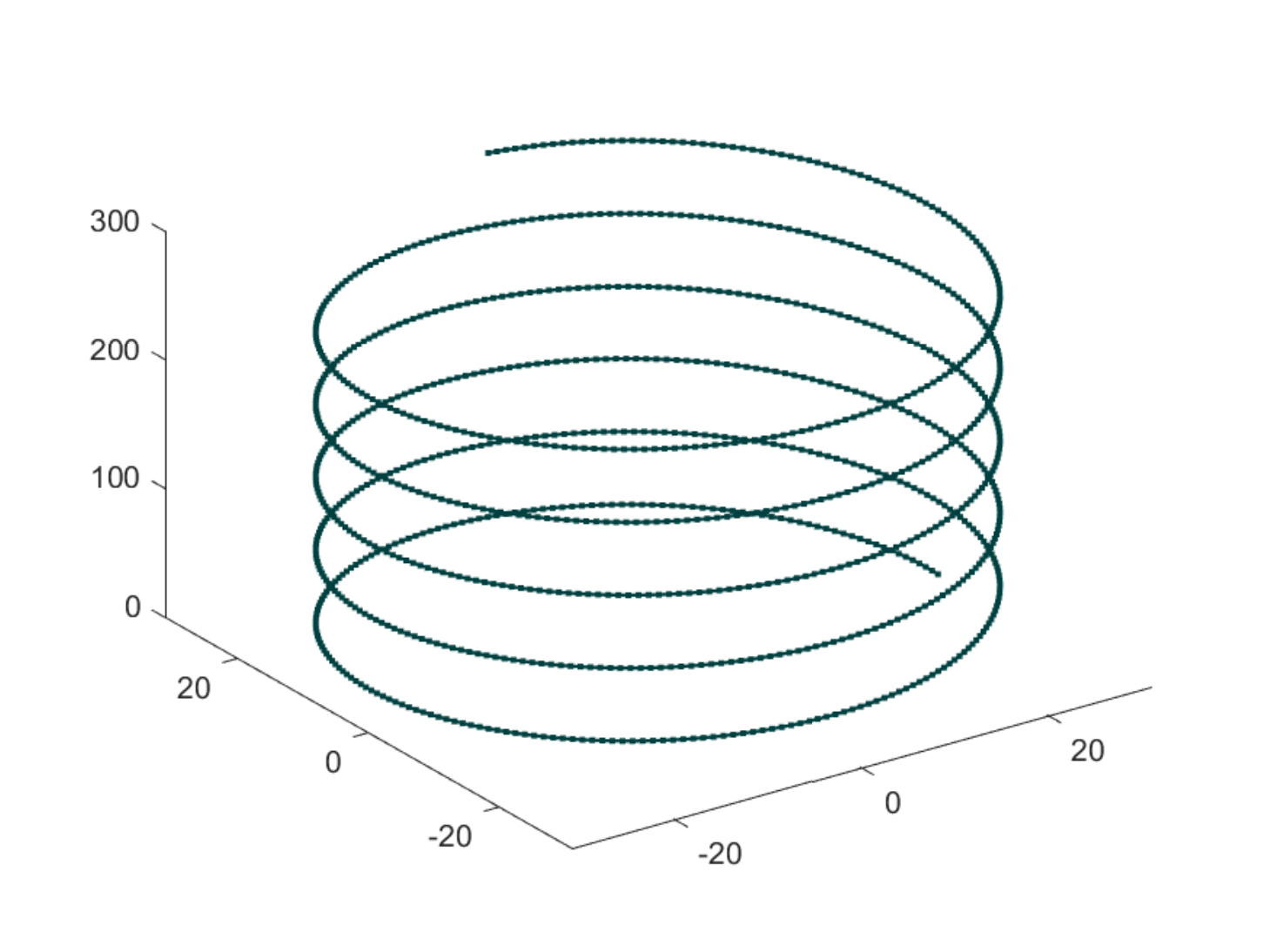}}
	\subfigure[The mFDBK curve]{
	    \includegraphics[width=0.48\textwidth]{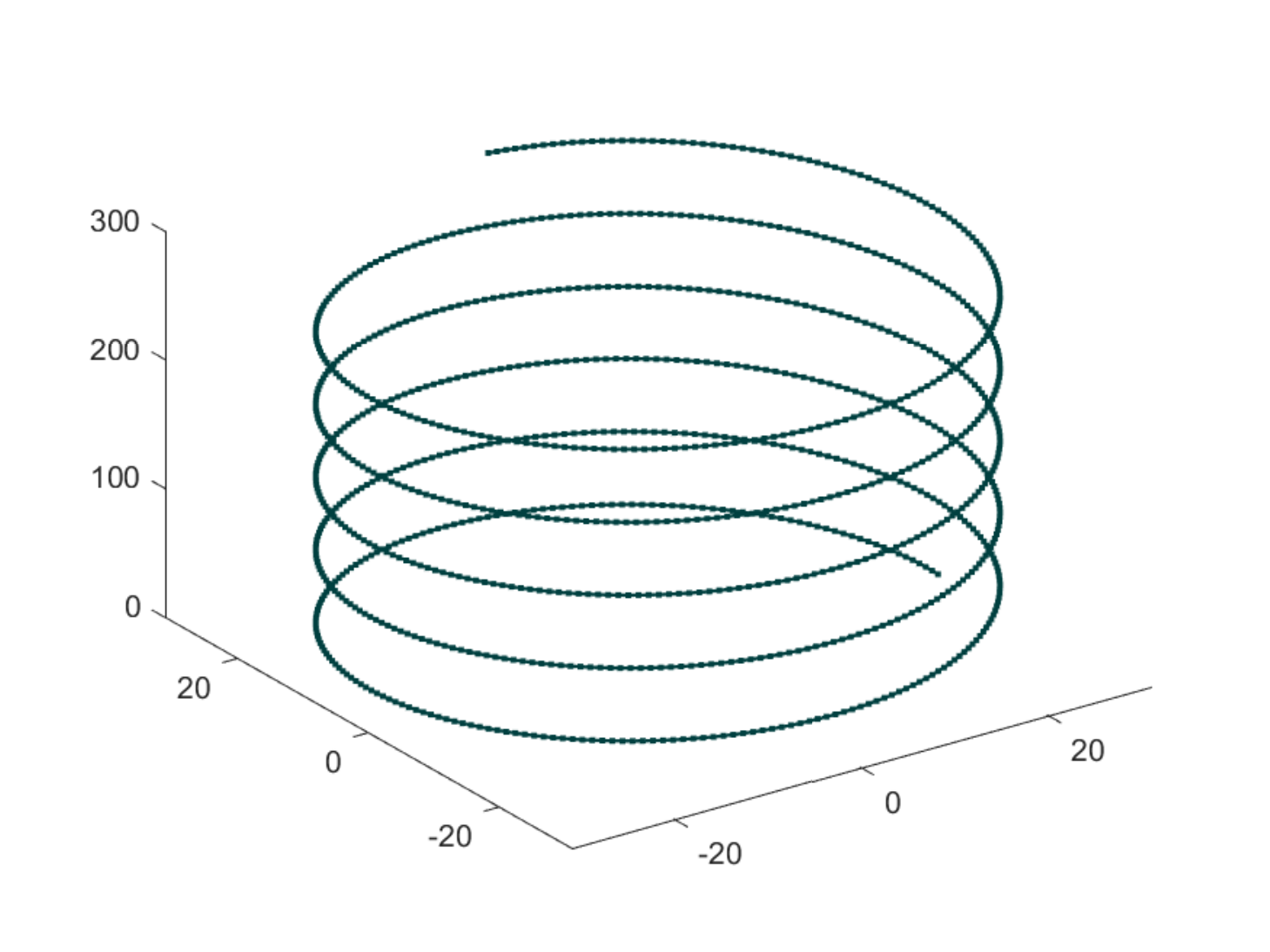}}
\caption{The limiting curves given by mMWRK (a) and mFDBK (b) for Example \ref{mMWRK and mFDBK:example2} when $m=20000$ and  $n=500$.}
\label{fig:FinalCurve1}
\end{figure}

\begin{figure}[!htb]
\centering
	\subfigure[The mMWRK curve]{
	    \includegraphics[width=0.48\textwidth]{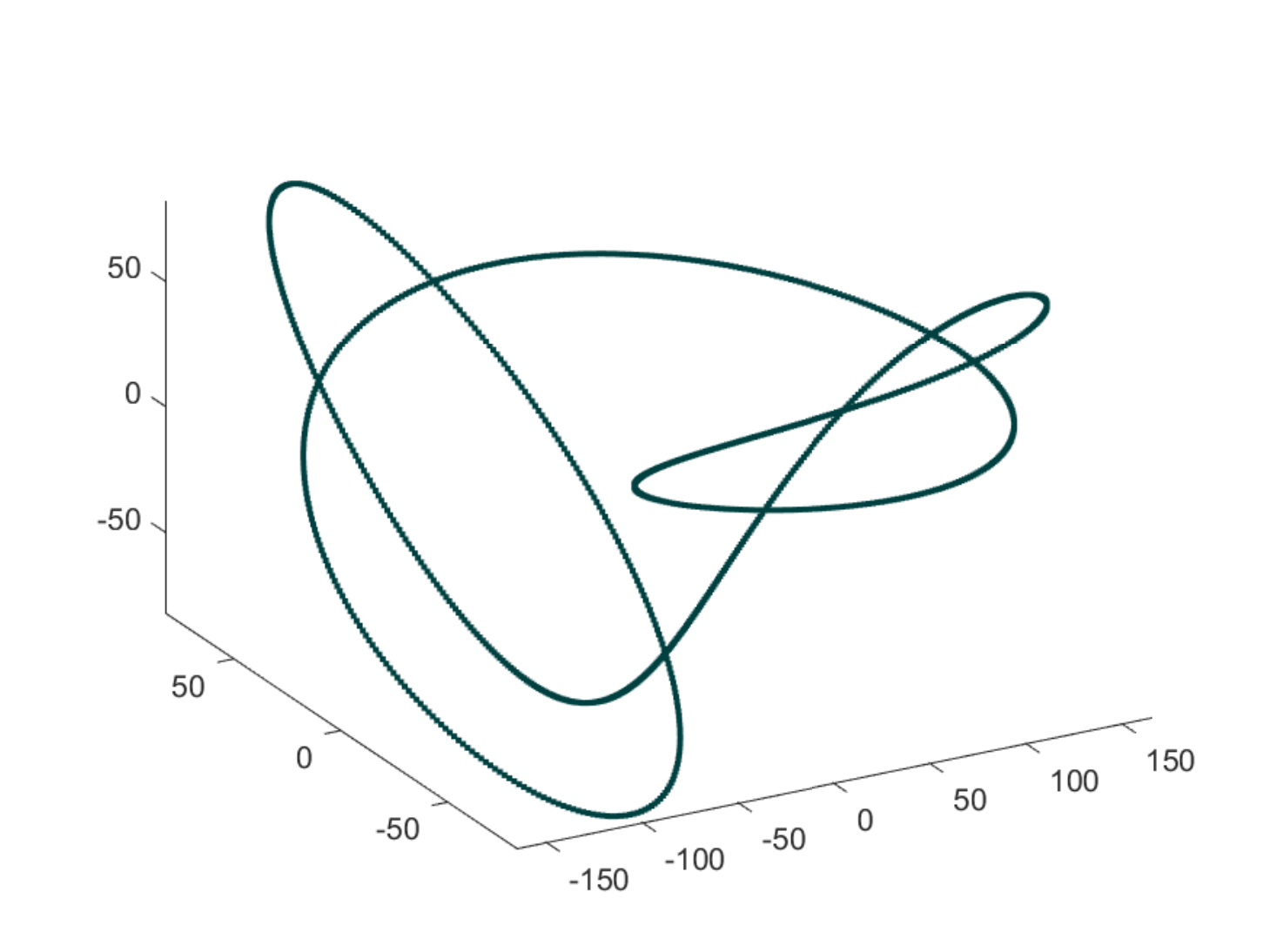}}
	\subfigure[The mFDBK curve]{
	    \includegraphics[width=0.48\textwidth]{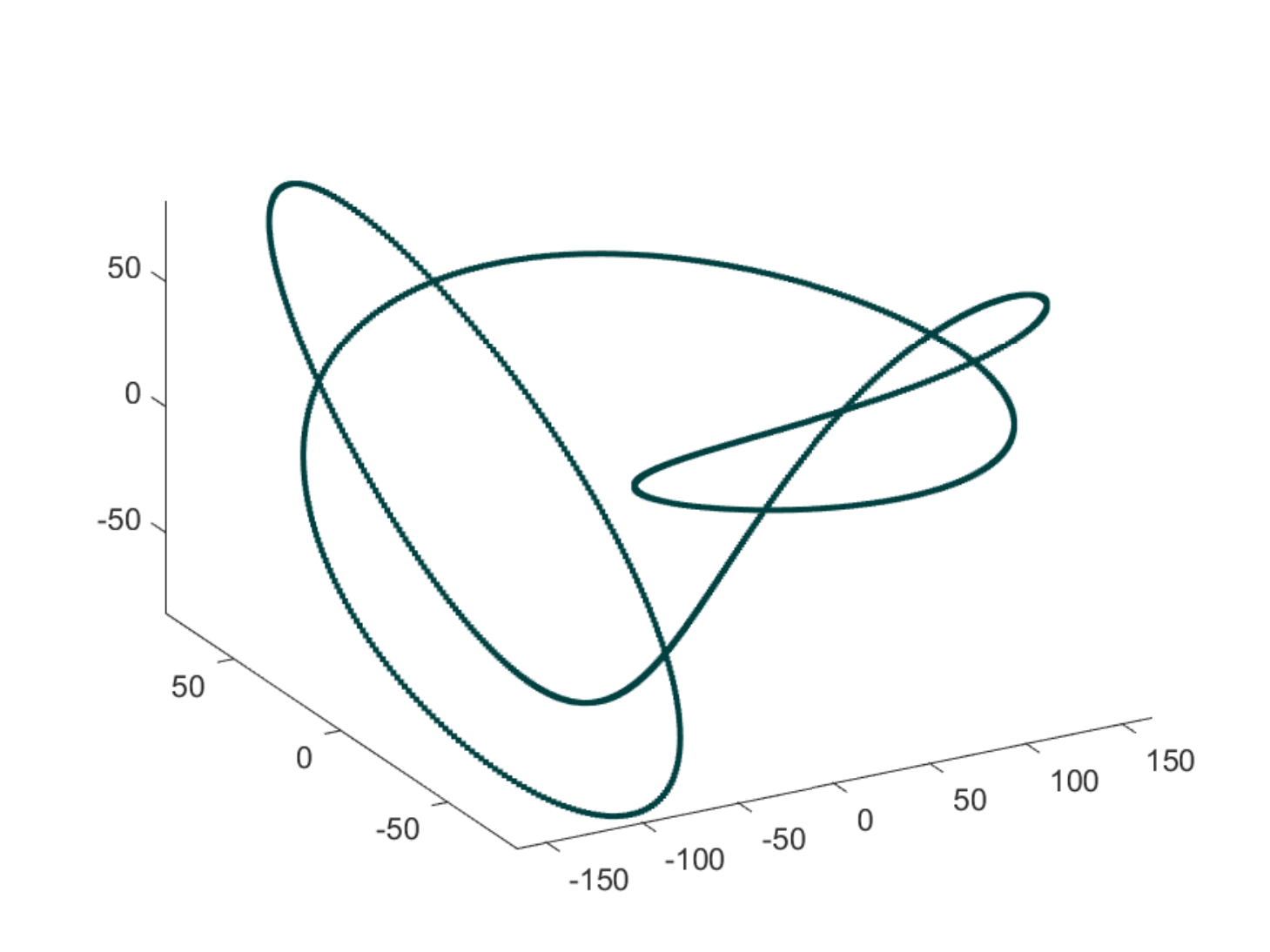}}
\caption{The limiting curves given by mMWRK (a) and mFDBK (b) for Example \ref{mMWRK and mFDBK:example2} when $m=20000$ and  $n=500$.}
\label{fig:FinalCurve2}
\end{figure}

We first list in Tables \ref{tab:Ex_curveFitting_Helical}-\ref{tab:Ex_curveFitting_GrannyKnot} the number of iteration steps and speed-ups of the mMWRK and mFDBK methods against their momentum-free variants. From these tables, we see that the mMWRK (resp., mFDBK) method outperforms
the MWRK (resp., FDBK) method. In particular, the speed-up is stable at $1.55$. We then plot the limiting curve of the mMWRK and mFDBK methods in Figures \ref{fig:FinalCurve1}-\ref{fig:FinalCurve2} for $m=20000$ and  $n=500$. As can be seen from the figures, both momentum methods achieve success in converging to the least-squares fitting curve.

\begin{table}[!htb]
\caption{The numerical results obtained by MWRK, mMWRK, FDBK, and mFDBK using $n$ control points to fit $m$ data points in Curve 1.}
\centering
\begin{tabular}{ccccc}
\cline{1-5}
$m\times n$&
$10000 \times 350$&$10000 \times 400$&
$10000 \times 450$&$10000 \times 500$\\
\cline{1-5}
MWRK  &5164&5810&6294&6646\\
mMWRK &3181&3664&4095&4657\\
SU$_1$&{\bf  1.62}&{\bf  1.59}&{\bf  1.54}&{\bf  1.43}\\
\cline{1-5}
FDBK  &844&853&870&921\\
mFDBK &510&537&582&612\\
SU$_2$&{\bf  1.65}&{\bf  1.59}&{\bf  1.49}&{\bf  1.50}\\
\cline{1-5}
$m\times n$&
$15000 \times 350$&$15000 \times 400$&
$15000 \times 450$&$15000 \times 500$\\
\cline{1-5}
MWRK  &5165&5823&6310&6811\\
mMWRK &3151&3654&4104&4636\\
SU$_1$&{\bf  1.64}&{\bf 1.59}&{\bf  1.54}&{\bf  1.47}\\
\cline{1-5}
FDBK  &807&812&916&916\\
mFDBK &528&547&584&610\\
SU$_2$&{\bf  1.53}&{\bf  1.48}&{\bf  1.57}&{\bf  1.50}\\
\cline{1-5}
$m\times n$&
$20000 \times 350$&$20000 \times 400$&
$20000 \times 450$&$20000 \times 500$\\
\cline{1-5}
MWRK  &5135&5762&6382&6751\\
mMWRK &3223&3620&4160&4587\\
SU$_1$&{\bf  1.59}&{\bf  1.59}&{\bf  1.53}&{\bf 1.47}\\
\cline{1-5}
FDBK  &778&893&828&954\\
mFDBK &537&543&567&629\\
SU$_2$&{\bf  1.45}&{\bf  1.64}&{\bf  1.46}&{\bf  1.52}\\
\cline{1-5}
\end{tabular}
\label{tab:Ex_curveFitting_Helical}
\end{table}

\begin{table}[!htb]
\caption{The numerical results obtained by MWRK, mMWRK, FDBK, and mFDBK using $n$ control points to fit $m$ data points in Curve 2.}
\centering
\begin{tabular}{ccccc}
\cline{1-5}
$m\times n$&
$10000 \times 350$&$10000 \times 400$&
$10000 \times 450$&$10000 \times 500$\\
\cline{1-5}
MWRK  &5250&6069&6352&7201\\
mMWRK &3284&3811&4235&4745\\
SU$_1$&{\bf 1.60 }&{\bf 1.59}&{\bf  1.50}&{\bf  1.52}\\
\cline{1-5}
FDBK  &838&907&960&954\\
mFDBK &572&597&622&673\\
SU$_2$&{\bf 1.47}&{\bf  1.52}&{\bf  1.54}&{\bf  1.42}\\
\cline{1-5}
$m\times n$&
$15000 \times 350$&$15000 \times 400$&
$15000 \times 450$&$15000 \times 500$\\
\cline{1-5}

MWRK  &5174&5932&6442&6976\\
mMWRK &3329&3792&4288&4764\\
SU$_1$&{\bf 1.55}&{\bf  1.56}&{\bf  1.50}&{\bf  1.46}\\
\cline{1-5}
FDBK  &850&925&989&904\\
mFDBK &570&599&583&642\\
SU$_2$&{\bf 1.49}&{\bf  1.54}&{\bf  1.70}&{\bf  1.41}\\
\cline{1-5}
$m\times n$&
$20000 \times 350$&$20000 \times 400$&
$20000 \times 450$&$20000 \times 500$\\
\cline{1-5}

MWRK  &5212&5947&6188&6959\\
mMWRK &3318&3812&4288&4747\\
SU$_1$&{\bf  1.57}&{\bf  1.56}&{\bf  1.44}&{\bf  1.47}\\
\cline{1-5}
FDBK  &832&925&923&982\\
mFDBK &572&601&595&636\\
SU$_2$&{\bf  1.45}&{\bf  1.54}&{\bf  1.55}&{\bf  1.54}\\
\cline{1-5}
\end{tabular}
\label{tab:Ex_curveFitting_GrannyKnot}
\end{table}

\section{Conclusions}\label{Sec:mMWRK+mFDBK+Conclusions}
For iteratively computing the minimum Euclidean-norm least squares solution of a consistent linear system, the MWRK and FDBK methods extend the deterministic single and multiple row-action methods, respectively, by introducing a promising adaptive index selection strategy. To further accelerate the convergence rate of the MWRK and FDBK methods, in this work, we utilize Polyak's heavy ball momentum acceleration technique and present the mMWRK and mFDBK methods and their computational complexity analysis. Convergence theory has been developed for the mMWRK and mFDBK methods. Some numerical examples, where the coefficient matrix is obtained from synthetic data and curve fitting, are given to demonstrate their numerical advantage over the MWRK and FDBK methods in terms of iteration counts. Numerical results illustrate that the mMWRK and mFDBK methods are competing row-action variants for solving the consistent linear systems.

Finally, we point out that finding the optimal values of step-size and momentum parameters in the mMWRK and mFDBK methods is a technical and skillful issue. These two parameters are determined by various factors such as the concrete structure and property of the coefficient matrix. This topic is of real value and theoretical importance. We will investigate this in detail in the future.

\section*{Acknowledgment}
This work is supported by the National Natural Science Foundation of China under grant 12201651.

\end{document}